\documentclass[11pt]{amsart}
\usepackage{latexsym}
\usepackage{amssymb,amscd,amsthm,amsfonts,verbatim,amsopn,color}
\usepackage[dvips]{epsfig}
\def\C{{\mathbb C}}

\def\R{{\mathbb R}}

\def\P{{\mathbb P}}
\def\Sphere{{\mathbb S}}

\def\BB{{\mathcal B}}

\def\FF{{\mathcal F}}
\def\GG{{\mathcal G}}
\def\HH{{\mathcal H}}
\def\LL{{\mathcal L}}

\def\NN{{\mathcal N}}

\def\Ss{{\mathcal S}}

\newcommand{\lra}{\longrightarrow}

\newcommand{\mfA}{\mathfrak A}

\newcommand{\lf}{\lfloor}
\newcommand{\rf}{\rfloor}

\newcommand{\wti}{\widetilde}

\newtheorem{thm}{Theorem}[section]
\newtheorem{df} [thm]{Definition}
\newtheorem{lm}  [thm]{Lemma}
\newtheorem{crl} [thm]{Corollary}
\newtheorem{prop}[thm]{Proposition}

\newtheorem{rem}[thm]{Remark}

\newtheorem{expl}[thm]{Example}
\newtheorem{alg}[thm]{Algorithm}
\numberwithin{equation}{section}

\newenvironment{explrm}{\begin{expl} \rm}{\end{expl}}

\newenvironment{dfrm}{\begin{df} \rm}{\end{df}}

\begin{document}

\title[Matroid polytopes, nested sets and  Bergman fans]
{Matroid polytopes, nested sets and  Bergman fans}

\author{Eva Maria Feichtner \, and \, Bernd Sturmfels}

\address{
Department of Mathematics, ETH Z\"urich, 8092 Z\"urich, Switzerland}
\email{feichtne@math.ethz.ch} 
\address{Department of Mathematics, UC Berkeley, Berkeley CA 94720, USA}
\email{bernd@math.berkeley.edu}

\date{\today}

\begin{abstract}
  The tropical variety defined by linear equations with constant
  coefficients is the Bergman fan of the corresponding matroid.
  Building on a self-contained introduction to matroid polytopes, we
  present a geometric construction of the Bergman fan, and we discuss
  its relationship with the simplicial complex of nested sets in the
  lattice of flats.  The Bergman complex is triangulated by the nested
  set complex, and the two complexes coincide if and only if every
  connected flat remains connected after contracting along any
  subflat.  This sharpens a result of Ardila-Klivans who showed that
  the Bergman complex is triangulated by the order complex of the
  lattice of flats.  The nested sets specify the De~Concini-Procesi
  compactification of the complement of a hyperplane arrangement,
  while the Bergman fan specifies the tropical compactification.
  These two compactifications are almost equal, and we highlight the
  subtle differences.
\end{abstract}


\maketitle


\section{Introduction}
\label{sec_intr}

Let $V$ be an $r$-dimensional vector subspace of the $n$-dimensional
vector space $\C^n$ over the field $\C$ of complex numbers. The {\em
  amoeba} of $V$ is the set of all vectors of the form
$$
\bigl( {\rm log} | v_1 |, {\rm log} | v_2 |, \ldots, {\rm log} |
v_n | \bigr) \,\,\in \,\, \R^n, $$
where $ (v_1,\ldots,v_n)$ runs over
all vectors in $V$ whose coordinates are non-zero. The asymptotic
behavior of the amoeba is given by an $r$-dimensional polyhedral fan
in $\R^n$. This fan was called the {\em logarithmic limit set of $V$}
in George Bergman's seminal paper~\cite{Be}. We use the term {\em
  Bergman fan} for this object. The study of such polyhedral spaces is
now an active area of research, known as {\em tropical geometry}
\cite{Hak, Spe, SS, Tev}. This paper is concerned with the {\em
  tropical variety} defined by a system of linear equations with
constant coefficients. The case of linear equations with coefficients
in a power series field is treated in \cite{Spe}. Bergman fans are the
local building blocks of Speyer's {\em tropical linear spaces}.

Our starting point is the observation of \cite[\S 9.3]{St1} that the
Bergman fan of a linear space $V$ depends only on the associated {\em
  matroid}. One way to specify this matroid is by its collection of
circuits $C$. These are the minimal sets arising as supports of linear
forms $\,\sum_{i \in C} a_i x_i\,$ which vanish on $V$. Introductory
references on matroids include \cite{BGW, Ox, Whi}.

  Let $M$ be any matroid of rank $r$ on 
the ground set $[n] = \{1,2,\ldots,n\}$. The {\em Bergman fan}
$\, \wti \BB(M)\,$ is the set of all vectors
$w = (w_1,\ldots,w_n) \in \R^n$ such that,
for every circuit~$C$ of $M$, the minimum of the set
$\,\bigl\{ \,w_i \,\, |\,\, i \in C \,\}\,$ is attained at least twice.
Note that $\wti \BB(M)$ is invariant under translation
along the line $\R(1,\ldots,1)$ in $\R^n$ and under positive
scaling. Hence we lose no information restricting our attention to
the intersection of the fan $\wti \BB(M)$ with the unit sphere 
in the hyperplane orthogonal to the line $\R(1,\ldots,1)$ in~$\R^n$:
\begin{equation}
\label{BSphere}
      \BB(M)\, \,
         := \, \,  \wti \BB(M)\,\cap\, \Sphere \, \qquad  \mbox{with }\, \, \,
      \Sphere \,=\,\bigl\{\,w \in \R^n\,|\,\sum_{i=1}^n\, w_i=0 \,\,\,\hbox{and} \,\,\,
      \sum_{i=1}^n \,w_i^2 =1\,\bigr\}\, . 
\end{equation}
The spherical set $\BB(M)$ is called the {\em Bergman complex\/} of
the matroid $M$. The terms ``fan'' and ``complex'' are justified by
our discussion of matroid polytopes in Section 2. In fact, $\wti
\BB(M)$ has a canonical decomposition as a subfan of the normal fan of
the matroid polytope, and, accordingly, $\BB(M)$ is a complex of
spherical polytopes.  Example \ref{cubearrange} reveals that the faces
of the Bergman complex are not always simplices.

In Section 3 we introduce the nested set complexes of an arbitrary lattice. 
If the lattice is the Boolean lattice $2^{[r]}$ of all subsets of
$[r]=\{1,2,\ldots,r\}$ then each nested set complex arises
as the boundary of a simplicial $(r-1)$-polytope.
The simple polytope dual to that simplicial polytope  is constructed
as the Minkowski sum of faces of the $(r-1)$-simplex;
we use this to develop a polyhedral theory of local Bergman 
complexes.

In Section 4 we apply the local theory of Section 3
to the lattice of flats $\LL_M$ of a rank $r$ matroid $M$,
and we  derive the theorem that
every nested set complex of $\LL_M$ is a unimodular 
triangulation of the Bergman complex $\BB(M)$.
We examine the local structure of this triangulation
in matroid-theoretic terms, thus refining the
results in \cite{AK}.

Among all nested set complexes of a matroid $M$,
there is always a minimal one, whose vertices
are indexed by the connected flats of $M$.
Section 5 is concerned with this minimal nested  set complex.
It is generally quite close to the Bergman complex $\BB(M)$,
and we show that they are equal if and only if every contraction
of a connected flat remains connected. We also discuss
algorithmic tools for computing the Bergman complex
along with its triangulation by minimal nested sets,
and we discuss some non-trivial examples.

In Section 6 we relate our combinatorial results to algebraic
geometry. The space $X\,{=}\, V{\cap}(\C^*)^n \,$ is the complement of
an arrangement of $n$ hyperplanes in complex affine $r$-space.  The
nested set complex specifies the {\em wonderful compactification} of
$X$, due to De~Concini and Procesi \cite{DP}, while the Bergman
complex specifies the {\em tropical compactification} of~$X$, due to
Tevelev \cite{Tev}.  The subdivision of Section 4 induces a canonical
morphism from the former onto the latter. We describe this morphism
geometrically.

  We close the introduction with two examples
  where our complexes are one-dimensional.

  \begin{expl} \label{K4} \rm
  $(n=6, r=3)$
  Let $M$ be the graphical matroid of the
  complete graph~$K_4$ on four vertices.
Here the nested set complex coincides with the
Bergman complex, and it equals the
{\em Petersen graph}, as depicted in Figure \ref{fig_NPi4Gmin} below.
When passing to the order complex of $\LL_M$,
three of the edges are subdivided into two, so
the order complex is a graph with $13$ vertices and $18$ edges.
For connections to phylogenetics see \cite[\S 3]{AK}.
  \end{expl}

  \begin{expl} \label{K4minus} \rm
  $(n=5, r=3)$ Consider the graph gotten from $K_4$ by removing
  one edge, and let $M'$ be the corresponding graphic matroid.
  The Bergman complex $\BB(M')$ is the {\em complete bipartite
  graph} $K_{3,3}$. Its six vertices are the two $3$-cycles
  and the four edges adjacent to the missing edge.
One of the nine edges connects the two vertices
indexed by the two $3$-cycles. In the nested set complex of $M'$
that edge is further subdivided by one vertex, corresponding
to the edge of $K_4$ which is disjoint from the missing edge.
This example appeared in \cite[Example 9.14]{St1} and
we shall return to it in Example \ref{ex_NPi4Gmin}.
   \end{expl}


\section{The matroid polytope and the Bergman complex}

We start with a brief introduction to matroid theory
with an emphasis on polyhedral aspects.
Let $M$ be a family of $r$-element subsets of
the ground set $[n] = \{1,2,\ldots,n\}$.  We represent each subset
$\sigma = \{\sigma_1, \ldots, \sigma_r \}$ by the 
corresponding sum of $r$ unit vectors
$$ e_\sigma \,\,\,\, = \,\,\,\, \sum_{i=1}^r e_{\sigma_i} \,\quad \in \,\,\, \R^n . $$
The set family $M$ is represented by the convex hull of these points
$$ P_M \quad := \quad {\rm conv}  \bigl\{ \,\,e_\sigma \,:\, \sigma \in M \,\,\bigr\}
\qquad \subset \, \,\, \R^n $$
This is a convex polytope of dimension $\leq n{-}1$. 
It is a subset of the $(n{-}1)$-simplex
$$ \Delta \quad = \quad \bigl\{ \,
(x_1,x_2,\ldots,x_n) \in \R^n \,\,: \,\, x_1 \geq 0, \ldots, x_n \geq 0, \,\,\,
x_1 + x_2 + \cdots + x_n = r \,\bigr\}. $$

\begin{df} \label{dfmatroid} A {\em matroid of rank $r$} is a family
$M$ of $r$-element subsets of $[n]$ such that
every edge of the polytope $P_M$ is
parallel to an edge of the simplex $\Delta$.
\end{df}

Experts in matroid theory may be surprised to see this unusual 
definition, but it is in fact equivalent to the many other definitions
familiar to combinatorialists. This equivalence was first proved by
Gel'fand, Goresky, MacPherson and Serganova \cite{GGMS}.
It forms the point of departure for the theory of
{\em Coxeter matroids}  in \cite{BGW}. This  suggests that
it would be worthwhile to extend the results in this paper to root systems
other than $A_n$.

The basic idea behind Definition \ref{dfmatroid}
  is as follows.  Every edge of the simplex $\Delta$
has the form ${\rm conv}(r e_i, r e_j)$, so it 
is parallel to a difference $e_i - e_j$  of two unit vectors.
The elements $\sigma \in M$ are the {\em bases} of the matroid,
and two bases $\sigma, \tau$ are connected by an edge
${\rm conv}(e_\sigma, e_\tau)$ if and only if
$e_\sigma - e_\tau = e_i - e_j$. The latter condition is equivalent to
$\,\sigma \backslash \tau = \{i\}$ and $\,\tau \backslash \sigma = \{j\}$,
so the edges of $P_M$ represent the basis exchange axiom.

Here is a brief summary of matroid terminology. Fix a matroid $M$ on $[n]$.
A subset $I \subseteq [n]$ is {\em independent} in $M$ if $I \subset \sigma$
for some basis $\sigma$. Otherwise $I$ is {\em dependent}.
The {\em rank} of a subset $F \subseteq [n]$ is the
cardinality of the largest independent subset of $F$.
A {\em circuit} is a dependent set which is minimal with respect to inclusion.
A subset $F \subseteq [n]$ is a {\em flat} of $M$ if there is no circuit $C$
such that $C \backslash F$ consists of exactly one element.
The intersection of two flats is again a flat. 
The span of $G \subseteq [n]$ is the smallest flat $F$ with $G \subseteq F$.
The collection of all flats is partially ordered by inclusion.
The resulting poset $\mathcal{L}_M$ is a {\em geometric lattice}, where 
$G_1 \wedge G_2 = G_1 \cap G_2$ and
$G_1 \vee G_2 = $ the span of $G_1 \cup G_2$.

The polytope $P_M$ is called the {\em matroid polytope} of $M$.
What we are interested in here is the following natural question concerning
the inclusion $\, P_M \subset \Delta$:
{\sl How does the boundary $\partial P_M$ of the matroid polytope 
intersect the boundary  $\partial \Delta$ of the ambient simplex~?}
The objects of study in this paper are polyhedral complexes
  which represent the topological space
$\,\partial P_M \backslash \partial \Delta$.
As we shall see, their combinatorial structure 
is truly ``wonderful''  \cite{DP}.

\begin{expl} \label{octahedron} \rm
Let $r = 2$ and $n = 4$ and consider the {\em uniform matroid} 
$ M = \bigl\{ \! \{1,2\}, \{1,3\}, \{1,4\}, \{2,3\}, \{2,4\}, \{3,4\} \! \bigr\} $.
Its matroid polytope is the regular octahedron
$$ P_M \,\,= \,\, {\rm conv}
\bigl\{\, (1,1,0,0), \, (1,0,1,0), \, (1,0,0,1), \, (0,1,1,0), \,
  (0,1,0,1), \, (0,0,1,1)\bigr\} \,\,\, \subset \,\,\, \R^4.  $$
  Here $\,\partial P_M \backslash \partial \Delta\,$
  consists of the relative interiors of  four
  of the eight triangles in $\partial P_M$. 
Here the Bergman complex consists of four points,
  which is the homotopy type of
  $\,\partial P_M \backslash \partial \Delta$.
  \end{expl}

In order to understand the combinatorics of
   $\,\partial P_M \backslash \partial \Delta\,$ for general matroids $M$,
we need to represent the matroid
polytope $P_M$ by a system of  linear inequalities.

\begin{prop}  \label{matroidineqs} The matroid polytope equals the following subset
of the simplex $\Delta$:
$$ P_M \quad = \quad \bigl\{ \, (x_1,\ldots,x_n) \in \Delta \,\,:\,\,
\sum_{i \in F} x_i \,\leq \, {\rm rank}(F) 
\quad \hbox{\rm for all flats} \,\,\, F \subseteq [n] \,\bigr\}. $$
\end{prop}

\begin{proof}
Consider any facet of the polytope $P_M$ and let
$\,\sum_{i=1}^n a_i x_i \leq b \,$ be an inequality defining this facet.
The normal vector $(a_1,a_2,\ldots,a_n)$  is perpendicular to
the edges of that facet.  But each edge of that facet is parallel
to some difference of unit vectors $e_i - e_j$. Hence
the only constraints on the coordinates of the normal vector are of 
the form $a_i = a_j$.
Using the equation $\,\sum_{i=1}^n x_i = r $ and scaling the
right hand side $b$, we can therefore assume that
$(a_1,a_2,\ldots,a_n)$ is a vector in $\{0,1\}^n$.
Hence the polytope $P_M$ is characterized by
the inequalities of the form $\,\sum_{i \in G} x_i \,\leq \, b_G \,$
for some $\,G \subseteq [n]$. The right hand side $b_G$ equals
$$ b_G \quad = \quad
{\rm max} \bigl\{ |\sigma| \cap G \,:\, \sigma \,\,\hbox{basis of}\,\, M \,\bigr\}
\quad = \quad {\rm rank}(G). $$
The second equality holds because every independent subset of $G$ can be completed
to a basis $\sigma$. Let $F$ be the flat spanned by $G$.  Then $G \subseteq F$
and ${\rm rank}(G) = {\rm rank}(F)$, and hence
the inequality $\,\sum_{ i \in F} x_i \leq {\rm rank}(F)\,$ implies
the inequality $\,\sum_{ i \in G} x_i \leq {\rm rank}(G) $.
\end{proof}

The circuit exchange axiom 
gives rise to the following equivalence relation on 
the ground set $[n]$ of the matroid $M$: $i$ and $j$ are {\em equivalent}
if there exists a circuit $C$ with $\{i,j\} \subseteq C$.
The equivalence classes are the {\em connected
components} of $M$. Let $c(M)$ denote the number
of connected components of $M$. We say that $M$
is connected if $c(M) =1 $.

\begin{prop} \label{dimformula}
The dimension of the matroid polytope $P_M$ equals $n-c(M)$.
\end{prop}

\begin{proof}
  Two elements $i$ and $j$ are equivalent if and only if there exist
  bases $\sigma$ and $\tau$ with $i\,{\in}\,\sigma$ and $\tau =
  (\sigma\backslash \{i\}) \cup \{j\}$. The linear space parallel to
  the affine span of $P_M$ is spanned by the vectors $e_i - e_j$
  arising in this manner. The dimension of this space equals~$n-c(M)$.
\end{proof}

Definition \ref{dfmatroid} implies that
every face of a matroid polytope is a matroid polytope.
Consider the face $P_{M_w}$ of $P_M$ at which 
the linear form $\,\sum_{i=1}^n w_i x_i\,$
attains its maximum. The bases of the matroid $M_w$
are precisely the bases  $\sigma$ of $M$ 
of maximal  $w$-cost $\, \sum_{i\in \sigma} w_i $.

Two vectors $w, w' \in \R^n$ are considered
{\em equivalent} for the matroid $M$ if $\,M_w = M_{w'}$.
The equivalence classes are relatively open convex
polyhedral cones. These cones form a complete
fan in $\R^n$. This fan is the {\em normal fan} of $P_M$.
If $\Gamma$ is a cone in the normal fan of $P_M$ and
$w \in \Gamma$ then we write $M_{\Gamma} = M_w$.
The following proposition shows that the Bergman fan
$\, \wti \BB(M)\,$ is a subfan of the normal fan of the matroid polytope $P_M$.

\begin{prop} \label{fiveprop}
The following are equivalent for a vector $w \in \R^n$:
\begin{enumerate}
\item The vector $w$ lies in the Bergman fan $\, \wti \BB(M)$.
\item The matroid $M_w$ has no loops.
\item Every element $i \in [n]$ lies in some basis of $M_w$.
\item The face $P_{M_w}$ has non-empty intersection with the interior
  of the simplex $\Delta$.
\item The linear functional $\,\sum_{i=1}^n w_i x_i\,$
attains its maximum over $P_M$ in $\, \partial P_M \backslash \partial \Delta$.
\end{enumerate}
\end{prop}

\begin{proof}
   The equivalence of (2) and (3) is the definition of loops in matroids
   namely, $i \in [n]$ is a {\em loop} of $M_w$ if it lies in no basis
   of $M_w$.  To see that (3) and (4) are equivalent, we note that the
   polytope $P_{M_w}$ is the convex hull of the vectors $e_\sigma$
   where $\sigma$ runs over all bases of $M_w$. This convex hull
   contains a point with full support (i.e.~intersects the interior of
   $\Delta)$ if and only if condition (3) holds. Conditions (4) and (5)
   are equivalent because $\,\sum_{i=1}^n w_i x_i\,$ attains its maximum
   over $P_M$ in $P_{M_w}$. Finally, condition (1) fails if and only if
   there exists a circuit $C$ of $M$ and an element $i \in C$ such that
   $w_i < {\rm min}\{ w_j \,: \,j \in C \backslash \{i\} \bigr\}$.
   This shows that if $i$ appears in any basis $\sigma$ of $M$ then we
   can replace $\sigma$ by another basis $ (\sigma \backslash \{i\})
   \,\cup \,\{j\} \,$ of higher $w$-weight, i.e., $i$ is a loop in
   $M_w$. Conversely, if (2) fails for $i$ then $M$ has a circuit $C$
   with $i \in C$ and $w_i < {\rm min}\{ w_j \,: \,j \in C \backslash
   \{i\} \bigr\}$, i.e.,~(1) fails.
  \end{proof}

  A question left open  in Proposition \ref{matroidineqs}
  is to identify the facet-defining inequalities for the matroid polytope $P_M$
other than the facets given by the ambient simplex $\Delta$.
  Let us assume that $M$ is connected, so
  that $P_M$ is $(n-1)$-dimensional.   We identify
  a flat $F$ with its $0/1$-incidence vector. Then
   $P_{M_F}$ is the face of $P_M$   at which the inequality
  $\,\sum_{i \in F} x_i \,\leq \, {\rm rank}(F) \,$ is attained
  with equality. A flat $F$ of $M$ is called a {\em flacet}
  if ${\rm dim}(P_{M_F}) = n-2$. Thus we can replace
  the word ``flat'' by the word ``flacet''
  in  Proposition \ref{matroidineqs}.

  In order to characterize the flacets of a matroid~$M$
  combinatorially, we note that every flat $F$ of rank~$s$
  defines two new matroids. The
  {\em restriction}  to $F$ is the rank~$s$ matroid
$$
  M [\emptyset,F] \quad = \quad \bigl\{ \, \sigma \,\cap \, F \,\,\, :\,\,\,
  \sigma \in M \,\,\hbox{and}\,\,|\sigma \,\cap \,F | = s \,\bigr\}.
  $$
  The {\em contraction} of $M$ at $F$ is the rank $r-s$ matroid
  $$ M[F, [n]] \quad = \quad \bigl\{\,
\sigma \backslash F \,\,\,: \,\,\, \sigma \in M \,\bigr\}. $$

\begin{prop} \label{flacetprop} A flat $F$ of a connected matroid $M$
is a flacet (i.e.~it defines a facet of the polytope $P_M$) if and only if 
both of the matroids $M[\emptyset,F]$ and $M[F,[n]]$ are connected.
\end{prop}

  \begin{proof} Let $s={\rm rank}(F)$.
  The matroid indexed by the incidence vector of the flat $F$ equals
\begin{equation}
  \label{bulletformula}
   M_F \quad = \quad
  M[\emptyset,F] \,\oplus \,M[F,[n]] \quad = \quad
\bigl\{ \,\sigma \in M \,\, : \,\,|\sigma \,\cap \,F | = s \,\bigr\}. 
\end{equation}
The corresponding face of $P_M$ is the direct product of matroid polytopes
$$ P_{M_F} \quad = \quad P_{M[\emptyset,F]} \,\times \,P_{M[F,[n]]} . $$
By Proposition \ref{dimformula}, the dimension of this polytope
equals
$$ n-c(M_F) \,\, = \,\, |F| - c(M[\emptyset,F]) \,\,+\,\, n-|F| - c(M[F,[n]]). $$
This number is at most $n-2$. It equals $n-2$ if and only if
$F$ is a flacet. This happens if and only if both  positive integers
$\,c(M[\emptyset,F]) \,$ and $\, c(M[F,[n]])\,$
are equal to one.
  \end{proof}

Assuming that $M$ is connected, the normal fan of $P_M$
defines a subdivision of  the $(n-2)$-sphere $\Sphere$ defined 
in (\ref{BSphere}).
This subdivision is isomorphic
to the polar dual $P_M^*$ of the matroid polytope.
(If $M$ is not connected then we replace
$\Sphere$ by an appropriate sphere of dimension
$n-c(M)-1$.)  Each cone $\Gamma$ in the normal fan of
$P_M$ is identified with a face of the 
{\em dual matroid polytope} $P_M^*$.
Hence every face $\Gamma$ of $P_M^*$ has  an 
associated rank~$r$ matroid $M_\Gamma$. We can 
now recast the definition of the Bergman complex as follows:

\begin{thm} \label{dfbergman} 
The Bergman complex $\BB(M)$ of a matroid $M$ is the
   subcomplex in the boundary $\partial P_M^*$ of the dual matroid
   polytope $P_M^*$ which consists of all faces $\Gamma$ such that
   $M_\Gamma$ has no loop.
     The vertices of $\BB(M)$ are indexed by the flacets of $M$,
     provided $M$ is connected.
    The Bergman complex $\BB(M)$ is homotopy
   equivalent to the space $\partial P_M \backslash \partial \Delta$.
   \end{thm}

\begin{proof}
Proposition (\ref{fiveprop}) implies  that a vector $w \in \partial P_M^*$
lies in a face of the Bergman complex $\BB(M)$ if and only if
$M_\Gamma$ has no loop. This proves the first assertion.
The second assertion follows from (the discussion prior to)
Proposition \ref{flacetprop}.
The third assertion follows from the general 
fact that a subcomplex in the boundary of a polytope
is homotopy equivalent to the corresponding co-complex in the
boundary of the dual polytope. This version of Alexander duality is proved 
using barycentric subdivisions.
\end{proof}

For instance, in Example \ref{octahedron}, the
dual matroid polytope $P_M^*$ is a regular $3$-cube,
and the Bergman complex  $\BB(M)$ consists of
four of the eight vertices of the cube $P_M^*$.
Here is an example which shows that
the Bergman complex $\BB(M)$ is generally not simplicial.

\begin{expl} 
\label{cubearrange}
\rm Let $M$ be the rank $4$ matroid on 
$[6] = \{1,2,\ldots,6\}$ whose non-bases are $A = 1234$, $B = 1356$ and
$C = 2456$.  This matroid is realized by the six vertices of a triangular
prism, or by the six vertices of an octahedron, or by the six 
facet planes of a regular $3$-cube. Figure~\ref{fig_cubearrgt}
shows a picture of this affine hyperplane arrangement.
The three non-bases correspond to the
three special  points
$A,B,C$ on the ``plane at infinity''.

\begin{figure}[ht]
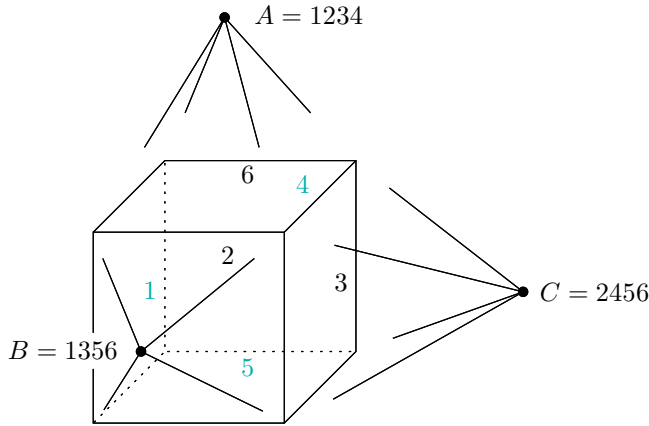

   \begin{picture}(0,0)%
     \includegraphics{cubearrgt.pstex}%
   \end{picture}%
   \input{cubearrgt.pstex_t}%
   
\caption{The cube realization of~$M$} 
\label{fig_cubearrgt}
\end{figure}

The matroid $M$ has six flacets of rank one, namely $1,2,3,4,5,6$,
and three flacets of rank two, namely $A,B,C$.
The Bergman complex $\BB(M)$ is two-dimensional. It has
$9$ vertices (the flacets),
$24$ edges, and $23$ two-dimensional faces.
Twenty faces are triangles:
$$
125 \quad
126\quad
145\quad
146\quad
235\quad
236\quad
345\quad
346 $$
$$
12A \,\,\, 14A \,\,\, 23A \,\,\, 34A \,\,\,
15B \,\,\, 16B \,\,\, 35B \,\,\, 36B \,\,\,
25C \,\,\, 26C \,\,\, 45C \,\,\, 46C $$
The remaining three two-dimensional faces are squares:
$$ 1A3B \qquad 2A4C \qquad 5B6C $$
The nested set complex, to be introduced in the next section,
will subdivide these squares by adding the diagonals
$13$, $24$ and $56$. We invite
the reader to compute the five-dimensional  matroid polytope  $P_M$
and to visualize the cocomplex $\partial P_M \backslash \partial \Delta$.
Here $\BB(M)$ and $\partial P_M \backslash \partial \Delta$ have the 
homotopy type
of a bouquet of seven two-dimensional spheres. 
\end{expl}


\section{Complexes of nested sets}\label{sect_ns}

Let $\LL$ be a finite lattice. We review the construction of a family
of simplicial complexes associated with $\LL$ due to Feichtner and
Kozlov~\cite{FK}, who generalized earlier work by De~Concini and
Procesi \cite{DP} on the special case (of interest here) when $\LL\,{
  =}\,\LL_M$ is the geometric lattice of a (realizable) matroid $M$.
Intervals in $\LL$ are denoted \mbox{$[X,Y]\,{=}\,\{Z\in \LL:$} $\,
X\,{\leq}\,Z\,{\leq}\,Y\}$. If $\,\Ss\,{\subseteq}\,\LL\,$ and
$\,X\,{\in}\,\LL\,$, then we write $\Ss_{\leq
  X}\,{:=}\,\{Y\,{\in}\,\Ss\,\,:\,\,Y\,{\leq}\,X\}$, and similarly for
$\Ss_{< X}$, $\Ss_{\geq X}$, and $\Ss_{> X}$ The set of maximal
elements in $\Ss$ is denoted $\mathrm{max}\, \Ss$.

\begin{dfrm} \label{df_building}
Let $\LL$ be a finite lattice.  A subset $\GG$ 
in~$\LL_{>\hat 0}$  is a {\em building set\/} if for any 
$X\,{\in}\,\LL_{>\hat 0}$  and
{\rm max}$\, \GG_{\leq X}=\{G_1,\ldots,G_k\}$ there is an isomorphism
of partially ordered sets
\begin{equation}\label{eq_buildg}
\varphi_X:\,\,\, \prod_{j=1}^k\,\,\, [\hat 0,G_j] \,\,
                                \stackrel{\cong}{\lra}
                                 \,\, [\hat 0,X]\, ,
\end{equation}
where the $j$-th component of the map $\varphi_X$ is the inclusion of 
intervals $\,[\hat 0,G_j]\,\subset \, [\hat 0,X]\,$ in $\LL$.
\end{dfrm}

The full lattice $\LL_{>\hat 0}$ is the simplest example of a building set 
for $\LL$. Besides this maximal building set 
(which we denote by $\LL$ for simplicity),
  there is always a minimal 
building set $\GG_{\mathrm{min}}$  consisting 
of all elements $X$ in $\LL_{>\hat 0}$ which do not allow for a 
product decomposition of the lower interval~$[\hat 0,X]$.
We call these elements the {\em connected elements\/} of~$\LL$. They 
have been termed {\em irreducible elements\/} at other places; we here
borrow the term {\em connected\/} from matroid terminology (see the beginning 
of Section~\ref{sect_NMvBM})

\begin{dfrm} \label{df_nested}
   Let $\LL$ be a finite lattice and $\GG$ a building set in $\LL$
   containing the maximal element $\hat 1$ of $\LL$.  A subset $\Ss$ in
   $\GG$ is called~{\em nested\/} if for any set of incomparable elements
   $X_1,\dots,X_t$ in $\Ss$ of cardinality at least two, the join
   $X_1\vee\dots\vee X_t$ does not belong to $\GG$.  The nested
   sets form an abstract simplicial complex $\wti \NN(\LL,\GG)$.
   Topologically, $\wti \NN(\LL,\GG)$ is a cone with apex $\hat 1$, its
   link  $\,\NN(\LL,\GG)$ is the {\em nested set complex\/} of
   $\LL$ with respect to $\GG$.
\end{dfrm}

If a building set $\GG$ does not contain the maximal element $\hat 1$,
it can always be extended to a building set $\widehat \GG\,{=}\,\GG
\,{\cup}\,\{\hat 1\}$. We can define nested sets with respect to~$\GG$
as above, and we then find the resulting abstract simplicial complex
being equal to the base of~$\wti \NN(\LL,\widehat \GG)$.

For the maximal building set $\LL$, the nested set
complex $\NN(\LL,\LL)$ coincides with the order complex $\Delta(\LL)$
of $\LL$, i.e., the abstract simplicial complex of totally ordered
subsets in the proper part of the lattice, $\,\LL\setminus\{\hat
0,\hat 1\}$. As we pointed out above, there is always the
unique minimal building set $\GG_{\mathrm{min}}$ of connected elements 
in~$\LL$, hence a simplicial
complex $\NN(\LL,\GG_{\mathrm{min}})$ with the least number of
vertices among the nested set complexes of~$\LL$.

\begin{explrm} \label{ex_NPi4Gmin}
Let $M$ be the matroid of the complete graph~$K_4$ as in 
Example~\ref{K4} and consider its lattice of flats $\LL_M$.
We depict its nested set complex with
   respect to the minimal building set $\GG_{\mathrm{min}}$ in
   Figure~\ref{fig_NPi4Gmin}.  The minimal building set is given by the
   atoms, the $3$-point lines  and the top
   element~$\hat 1$ in~$\LL_{M}$.
   Larger building sets are obtained by adding 
some of the three $2$-point lines.
This results in a subdivision of edges
in $\NN(\LL_M,\GG_{\mathrm{min}})$.

\begin{figure}[ht]
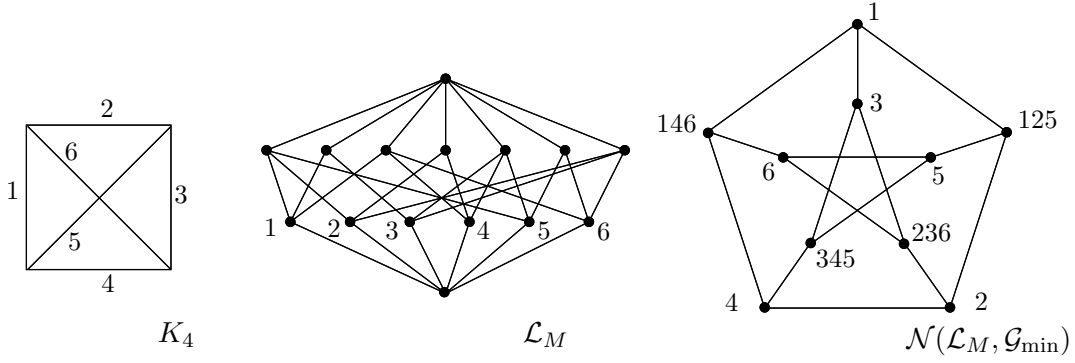

   \begin{picture}(0,0)%
     \includegraphics{BgmMK4.pstex}%
   \end{picture}%
   \input{BgmMK4.pstex_t}%
   
\caption{The lattice of flats $\LL_{M}$ and its minimal nested set 
complex $\NN(\LL_{M},\GG_{\mathrm{min}})$.} 
\label{fig_NPi4Gmin}
\end{figure}
\end{explrm}

\begin{explrm} \label{ex_NK4minus}
Remove one edge from the graph $K_4$
to get the graphic matroid $M'$ of Example~\ref{K4minus}.
  We depict its lattice of flats 
$\LL_{M'}$ in Figure~\ref{fig_NK4minus}. Again, the minimal 
building set $\GG_{\rm min}$ is given by the atoms, the two $3$-point lines 
and the top element. Any other building set is obtained by 
adding some of the four $2$-point lines. We depict the nested
set complex with respect to 
$\GG_{\rm min}$: it is a $K_{3,3}$ which is subdivided in 
one edge. Nested set complexes for larger building sets are
obtained by subdividing up to four further edges.

\begin{figure}[ht]
   \begin{picture}(0,0)%
     \includegraphics{NK4minus.pstex}%
   \end{picture}%
   \input{NK4minus.pstex_t}%
   
\caption{The lattice of flats $\LL_{M'}$ and its nested set 
complex $\NN(\LL_{M'},\GG_{\rm min})$.} 
\label{fig_NK4minus}
\end{figure}
\end{explrm}

\smallskip

A lattice $\LL$ is {\em atomic} if every element is a join of atoms.
The lattice of flats $\LL_M$ of a matroid is an
  atomic lattice. For arbitrary atomic lattices $\LL$, 
Feichtner and Yuzvinsky~\cite{FY} proposed the following
polyhedral realization of the nested set complexes 
$\NN(\LL,\GG)$.

\begin{dfrm} \label{df_nestedfan}
Let $\LL$ be an atomic lattice and $\mfA\,{=}\,\{A_1,\ldots,A_n\}$ its 
set of atoms; $\GG$ a building set in $\LL$ containing $\hat 1$. 
For any $G\,{\in}\,\GG$, let 
$\lf G \rf\,{:=}\,\{A\in \mfA\,|\, A\leq G\}$, the subset of atoms below~$G$.
We define the characteristic vector $v_G$ in~$\R^n$
by
\[
    (v_G)_i\, \, :=\, \, \left\{
\begin{array} {ll}
1 & \mbox{ if }\, A_i\in \,\lf G \rf  \\
0 & \mbox{ otherwise\,  }
\end{array}
\right.   \qquad \mbox{ for }\, i=1,\ldots, n\, .
\] 
For any nested set $\,\Ss\in \wti \NN(\LL,\GG) $, the set of vectors 
$\, \{v_G\,|\, G\in \Ss\}\,$ is  linearly independent, and it hence
spans a simplicial cone $\, V_{\Ss} \, = \,
\R_{\geq 0} \{v_G\,|\, G\in \Ss\}$.
These cones intersect along faces, namely
$\, V_{\Ss}\,\cap \,V_{\Ss'}\,= \,
V_{\Ss \,\cap  \,\Ss'}$, and hence they form
a simplicial fan $\,\wti\Sigma(\LL,\GG)\,$ in $\R^n$.

As we did before with the facet normals to the matroid polytope $P_M$,
  we can replace the vectors $v_G$ by equivalent vectors 
which lie on the
$(n{-}2)$-dimensional sphere
$\,\Sphere\,{=}\,\{\,w \,{\in}\, \R^n\,|\,
\sum_{i=1}^n w_i=0,\, \sum_{i=1}^n w_i^2 = 1\}$.
This is accomplished by translating $v_G$ along the
line $\R(1,\ldots,1)$ and then scaling it to have unit length.
We lose no information by replacing
$\,\wti \Sigma(\LL,\GG)\,$ with its restriction to the $(n{-}2)$-sphere 
$\Sphere$.
The resulting complex $\,\Sigma(\LL,\GG)\,$ is
a geometric realization of the 
abstract simplicial complex~$\NN(\LL,\GG)$.
\end{dfrm}

We recall some results 
concerning the geometry and topology of nested set complexes:

\begin{prop} \label{prop_subd} 
{\rm (\cite[Prop.\ 2]{FY}, \cite[Thm.\ 4.2, Cor.\ 4.3]{FM})}\newline
{\bf(1)}
For any atomic lattice $\LL$ and any building set $\GG$ in~$\LL$, the fan 
$\,\wti \Sigma(\LL,\GG)\,$ is unimodular. \newline
{\bf(2)}
For building sets $\HH\,{\subseteq}\,\GG$ in $\LL$,
the simplicial fan $\,\wti \Sigma(\LL,\GG)\,$ can be obtained from 
$\,\wti \Sigma(\LL,\HH)\,$ by a sequence of stellar subdivisions. 
In particular, the support sets of the fans $\,\wti \Sigma(\LL,\GG)\,$ 
coincide for all building sets $\GG$ in $\LL$. \newline
{\bf(3)} For an atomic lattice $\LL$ and any building set $\GG$ in $\LL$,
the nested set complex $\NN(\LL,\GG)$ is homeomorphic to the order complex 
$\Delta(\LL)$.
\end{prop}

We shall see in Theorem \ref{Main} that, for any 
matroid $M$ and any building set
$\GG$ in $\LL_M$, the fan $\,\wti \Sigma(\LL_M,\GG)\,$ is a 
refinement of the Bergman fan~$\,\wti \BB(M)\,$. In particular, 
the fans have the same support sets in~$\R^n$.
Also, the nested set complex $\, \NN(\LL_M,\GG)\,$ is
a triangulation of the Bergman complex $\BB(M)$ for any building 
set $\GG$ in $\LL_M$.
The special case of this result when $\GG\,{=}\,\LL_M$ 
and the nested set complex equals the order complex
of~$\LL_M$ is due to Ardila and Klivans \cite[Sect.\ 2, Thm.\ 1]{AK}.
It is possible to derive Theorem~\ref{Main} from their result
using the techniques of combinatorial blow-ups developed in~\cite{FK}.
However, we have chosen a different route which 
will keep this paper  self-contained.

\medskip

In what follows we take $\LL$ to be the Boolean lattice $2^{[r]}$
whose elements are the subsets of $[r] = \{1,2,\ldots,r\}$. Clearly,
$\,2^{[r]}\,$ is an atomic lattice, in fact, it is the lattice of
flats of the free rank $r$ matroid $\, M = \bigl\{ \{1,2,\ldots,r\}
\bigr\}$.  We will show (in Theorem \ref{einsfuenf}) that, for any
building set $\GG$ in~$2^{[r]}$, the fan
$\,\widetilde \Sigma(2^{[r]},\GG)$, when regarded modulo the line $\R(1,\ldots,1)$
as always, is the normal fan to a
simple $(r-1)$-polytope $\Delta_\GG$. Equivalently, $\, \NN(2^{[r]},\GG)\,$ is the
boundary of a simplicial $(r-1)$-polytope $\Delta_\GG^*$.  This dual
pair of polytopes should be of independent interest for further study
of the toric manifolds introduced in \cite{FY}.

\begin{rem} \rm
The minimal building set in the Boolean lattice
$2^{[r]}$ is the set of atoms, and, following our convention when defining 
nested set complexes, we include the maximal 
lattice element~$[r]$; hence $\,\GG_{\mathrm{min}} \, = \,
\bigl\{ \{1\}, \{2\}, \ldots,\{r\}, [r] \bigr\}$.
The fan $\wti\Sigma(2^{[r]},\GG_{\mathrm{min}})$ is the normal fan
to the $(r-1)$-simplex $\Delta_{\GG_{\mathrm{min}}}$, and
$\, \NN(2^{[r]},{\GG_{\mathrm{min}}})\,$ is the boundary complex
of the dual simplex $\,\Delta_{\GG_{\mathrm{min}}}^*$. On
the other extreme, the maximal nested set complex
$\NN(2^{[r]},\GG_{\mathrm{max}})$, $\,\GG_{\mathrm{max}}\,{=}\,2^{[r]}\,$, 
is the barycentric subdivision 
of the boundary of the $(r-1)$-simplex.
The corresponding fan  $\,\wti \Sigma(2^{[r]},\GG_{\mathrm{max}})\,$ 
is the braid arrangement $\{x_i = x_j\}$, and the
simple polytope $\Delta_{\GG_{\mathrm{max}}}$ is the
{\em permutohedron}. Hence the polytopes
$\Delta_\GG$  for $\GG$ building sets in $2^{[r]}$ interpolate between the
$(r-1)$-simplex $\Delta_{\GG_{\mathrm{min}}}$ and the
$(r-1)$-dimensional  permutohedron $\Delta_{\GG_{\mathrm{max}}}$.
This class of polytopes includes many interesting polytopes,
such as the {\em associahedron}, where $\GG$ is the set of all 
segments $\{i,i+1,\ldots,j-1,j\}$ for $1 \leq i < j \leq r $, and
the {\em cyclohedron}, where $\GG$ is the set of all cyclic segments.
Both of these polytopes are special cases of the
{\em graph-associahedra} of
Carr and Devadoss \cite{CD}.
\end{rem}

\begin{rem} \rm
After completion of this paper, we learned that most
of the results in the remainder of Section 3 had been
obtained independently by A.~Postnikov \cite{Pos}.
The main focus of Postnikov's work is the study of
the Ehrhart polynomials
of the polytopes $\Delta_\GG$.
\end{rem}

The following lemma characterizes arbitrary  building sets in
a Boolean lattice.

\begin{lm} \label{isbuilding}
A family $\FF$ of subsets of $[r]$ is
a building set in the Boolean lattice  $2^{[r]}$
  if and only if $\FF$ contains all singletons $\{i\}$, $i\in [r]$, and 
the following condition holds:
if $F,F' \in \FF$ and
$\,F \,\cap \, F' \not= \emptyset\,$
then $\,F \, \cup \, F' \,\in\, \FF$.
\end{lm}

\begin{proof}
If $\FF$ is a subset of $2^{[r]}$ and $X \in 2^{[r]}$ then $\FF_{\leq
   X}$ consists of all subsets of $X$ that are in the family $\FF$. The
poset map in (\ref{eq_buildg}) is an isomorphism if the factors
$G_1,\ldots, G_k$ are pairwise disjoint and their union is~$X$. We
conclude that $\FF$ is a building set for $2^{[r]}$ if and only if all
singletons are in~$\FF$ and, for every $X \in 2^{[r]}$, the
maximal elements in $\, \FF_{\leq X}\,$ are pairwise disjoint. This
condition is equivalent to the one stated in the lemma.
\end{proof}

\begin{lm}
Any family $\FF$ of subsets of $[r]$ can be enlarged
to a unique minimal family $\widehat \FF$ such 
that $\widehat \FF$ is a  building set in $2^{[r]}$.
We call $\widehat \FF$ the {\em building closure} of $\FF$.
\end{lm}

\begin{proof}
Fix a subset $X$ of $[r]$. We regard the set family
$\,{\rm max}\, \FF_{\leq X}\,$ as (the set of facets of) a 
simplicial complex of $[r]$. Let $\,\widehat \FF\,$
be the set of all subsets $X$ of $[r]$ such that
$X$ is a singleton or  the
simplicial complex $\,{\rm max}\, \FF_{\leq X}\,$ is connected.
It follows from Lemma \ref{isbuilding} that $\widehat \FF$
is a building set, and  every other building
set containing $\FF$ also contains $\widehat \FF$.
\end{proof}

We consider the standard simplex of dimension $r-1$. It is here denoted
$$ \Delta_{[r]} \quad =  \quad
\bigl\{ (x_1,\ldots,x_r) \in \R^r \,\,:\,\,
{\rm all} \,\,\,x_i \geq 0  \,\,\,
\hbox{and} \,\,\,
x_1 + x_2 + \cdots + x_r  = 1 \,\bigr\}. $$
Every subset $F$ of $[r]$ defines a face of $\Delta_{[r]}$ 
which is  a simplex of dimension $|F|-1$:
$$ \Delta_F \quad = \quad 
\bigl\{ (x_1,\ldots,x_r) \in \Delta_{[r]} \,\,:\,\,
x_i = 0 \,\,\hbox{for} \,\, i \not\in F \,\bigr\}. $$
With a family $\FF$ of subsets of $[r]$ we associate
the following Minkowski sum of simplices
$$\Delta_\FF \quad = \quad \sum_{F \in \FF} \Delta_F . $$
The dimension of the convex polytope
  $\Delta_\FF$ is given by the following formula:

\begin{rem} \label{dimpoly}
The dimension of the polytope $\Delta_\FF$ equals
$r-c$, where $c$ is the number of connected components
of the simplicial complex with facets $\,{\rm max} \,\FF$. 
\end{rem}

Each edge of the polytope $\Delta_\FF$ is parallel
to the difference of unit vectors $e_i - e_j$ in $\R^r$.
This means that $\Delta_\FF$ is a Minkowski
summand of the $(r-1)$-dimensional permutohedron.
In fact, $\Delta_\FF$ is a permutohedron whenever
$\FF$ contains all the two-element subset of $[r]$.
This observation implies the following facet description 
of $\Delta_\FF$.

\begin{prop} \label{ineqrepp} The polytope $\Delta_\FF$
consists of all non-negative vectors $(x_1,\ldots,x_r)$ such that
$\,x_1 + \cdots + x_r = | \FF| \,$ and the following inequality 
holds for all  subsets $\,G \,$ of $[r]$:
\begin{equation}
\label{facetineq}
  \sum_{i \in G} x_i \quad \geq \quad | \{ \, F \in \FF \, : \, F \subseteq G \,\bigr\}| .
  \end{equation}
Here it suffices to take those subsets $G$ which lie in the
  building closure $\,\widehat \FF\,$ of $\FF$.
\end{prop}

\begin{proof}
Since $\Delta_\FF$ is a Minkowski summand of the permutohedron,
it is defined by inequalities of the form
$\, \sum_{i \in G} x_i \,\,\geq \,\, \delta_G \,$ for some
parameters $\delta_G$. The minimum value of the
linear form $\sum_{i \in G} x_i$ on the simplex equals
one if $F \subseteq G$, and it equals zero otherwise.
This shows that $\,\delta_G =  | \{ \, F \in \FF \, : \, F \subseteq G \,\bigr\}|\,$
as desired. The face of $\Delta_\FF$ at which the linear form
$\, \sum_{i \in G} x_i \,$ attains its minimum can be expressed as follows:
\begin{equation}
\label{facett}
  \Delta_{\{\, F \in \FF \, : \, F \subseteq G \,\}} \quad + \quad
\Delta_{\{ F \backslash G \,: \,F \in \FF, \, F \backslash G \not= \emptyset \} } . 
\end{equation}
In order for this face to have codimension one in  $\Delta_\FF$, it is
necessary, by Remark \ref{dimpoly}, that the set family $\{F \in \FF \, : \, F \subseteq G \}$ represents
a connected simplicial complex. This condition is equivalent to 
the set $\,G \,$ being in the building closure $\, \widehat \FF$.
\end{proof}

\begin{crl} \label{crlirr}
If $[r] \in \FF$ then $\Delta_\FF$ is $(r-1)$-dimensional and its facets 
are indexed by the building closure $\widehat \FF$, i.e., the
inequality presentation in Proposition \ref{ineqrepp}
is irredundant.
  \end{crl}

\begin{proof}
We have ${\rm dim}(\Delta_\FF) = r-1$ by Remark \ref{dimpoly}.
Let $G \in \widehat \FF$. The left polytope in
  (\ref{facett}) has dimension $|G|-1$ as argued above.
  The right polytope in (\ref{facett}) contains $\Delta_{[r]\backslash G}$
  as a Minkowski summand, so it has dimension $r - |G|-1$. Hence the dimension
  of the face $(\ref{facett})$ is $\, (|G|-1)  + (r-|G|-1) \, = \, r-2$, 
which means it is a facet. 
\end{proof}

We are interested in conditions under which the
polytope $\Delta_{\FF}$ is simple, or equivalently,
the normal fan of $\Delta_{\FF}$ is simplicial.
The next theorem says that this happens if $\FF = \widehat \FF$.

\begin{thm} \label{einsfuenf} Let $\FF$ be a building set
in the Boolean lattice $2^{[r]}$ such
that $[r] \in \FF$. Then
$\Delta_\FF$ is an $(r-1)$-dimensional simple polytope,
and its normal fan is a unimodular simplicial fan which is
combinatorially isomorphic to the nested set complex
$\NN(2^{[r]}, \FF)$.
\end{thm}

\begin{proof}
Our assumptions say that $\widehat \FF = \FF$ and $[r] \in \FF$.
By Corollary \ref{crlirr}, the polytope~$\Delta_\FF$ has dimension $r-1$.
The facets of $\Delta_\FF$ are indexed
by the elements $G $ of $ \FF \backslash \{[r]\} $. The facet-defining inequality indexed by $G$
is  given in (\ref{facetineq}). We need to show that precisely
$r-1$ of these inequalities are attained with equality at any vertex of $\Delta_\FF$.

Pick a generic vector $w = (w_1,\ldots,w_r)$ and let  $v = (v_1,\ldots,v_r)$
be the vertex of~$\Delta_\FF$ at which $\sum_{i=1}^r w_i x_i$
attains its minimum. After relabeling, we may assume that
$w_1 < w_2 < \ldots < w_r$. Then the $i$-th coordinate of the vertex $v$ equals
$$ v_i \quad = \quad |\, \{\, F \in \FF \,\,: \,\, {\rm min}(F) \, = \, i \,\}\,|. $$
The inequality (\ref{facetineq}) indexed by $G \in \FF$
holds with equality at $v$ if and only if
$$ |\, \{\, F \in \FF \,\,: \,\, {\rm min}(F) \, \in \, G \,\}\,| \quad = \quad
\,\, | \{ \, F \in \FF \, : \, F \subseteq G \,\bigr\}|  $$
Given that $\FF$ is a building set, a
necessary and sufficient condition for this equality to hold is that the
set $G$ has the following specific form for some index $i \in [r]$:
$$ G_i \quad := \quad \bigcup \{ F \in \FF \, : \,\, {\rm min}(F) = i \,\} .$$
Here we are using the fact that $\FF$ is a building set,
which ensures that $G_i$ is in $\FF$.
Since $[r] \in \FF$, we have $G_1 = [r]$, which is excluded from
the sets in (\ref{facetineq}). Hence the facets incident to
$v$ are precisely the facets defined by $G_2,G_3,\ldots,G_r$.
In particular, $v$ is a simple vertex, and, by the relabeling argument
we conclude that $\Delta_\FF$ is a simple polytope.

The family $\{G_2,\ldots,G_r\}$ is a simplex in the
nested set complex $\NN(2^{[r]}, \FF)$. It remains to be
seen that this simplex is maximal and that all maximal 
simplices arise in this manner, after some permutation of $[r]$.
Indeed, suppose that $\Ss \subset \FF$ is a
facet of $\NN(2^{[r]}, \FF)$. 
The maximal elements of $\Ss$ are pairwise disjoint,
and since $[r] \in \FF$, their union has
cardinality less than $r$. After relabeling
we may assume that the element $1$ is not
in the union.  Since~$\Ss$ was assumed
to be maximal, its union equals
$\{2,\ldots,r\}$.  After relabeling again,
we can write the maximal sets as
$G_2, G_3,\ldots,G_k$ with  ${\rm min}(G_i) = i$.
Since the $G_i$ are pairwise disjoint, we can now
apply this construction recursively by restricting~$\Ss$ 
and~$\FF$ to each of the subsets $G_i$.
This  construction shows that
$\Ss$ has cardinality $r-1$, and it arises
precisely in the manner indicated above.
(See also Proposition  \ref{treerepresentation}
below).

It remains to note that the simple polytope
$\Delta_\FF$ is ``smooth'' in the sense of
toric geometry. The $r-1$ edges emanating
from the vertex $v$ have directions $e_i-e_j$,
and since the configuration of all vectors
$e_i-e_j$ is unimodular, it follows that
these $r-1$ edges form a basis for
the ambient lattice. It follows that the
normal fan of $\Delta_\FF$ is unimodular.
\end{proof}

Consider now an arbitrary  family
$\FF$ of subsets of $[r]$ and let
$\widehat \FF$ be  its building closure as before.
Since $\,\FF\subseteq \widehat \FF $, the polytope
$\Delta_\FF$  is a Minkowski summand of the
simple polytope~$\Delta_{\widehat \FF}$. 
The Minkowski summand relation of
convex polytopes corresponds to refinement
at the level of normal fans. Hence
Theorem \ref{einsfuenf} implies

\begin{crl} \label{normaltriang}
The normal fan of $\Delta_{\widehat \FF}$
is a triangulation of the normal fan of 
$\Delta_\FF$.
\end{crl}

\begin{explrm}
Let $r =4$ and $\FF = \bigl\{ 
\{1,2\},
\{2,3\},
\{3,4\},
\{1,4\} \bigr\}$. Then
$$ {\widehat \FF} \quad = \quad \FF \,\cup \,
\bigl\{ \{1\}, \{2\},\{3\}, \{4\}, \{1,2,3\},
\{1,2,4\},
\{1,3,4\},
\{2,3,4\},
\{1,2,3,4\} \bigr\}, $$
We consider orbits with respect to the
action of the dihedral group $D_4$  on~$\FF$ and 
on~$\widehat \FF$. The polytope
$\Delta_\FF$ is a three-dimensional
zonotope with  four zones, known as the
{\em rhombic dodecahedron}. It has
$14$ vertices, namely eight simple
vertices like  $(0,1,1,2)$, 
four four-valent vertices
like $(0,1,2,1)$, and two
four-valent vertices like  $(0,2,0,2)$.
The normal fan of $\Delta_\FF$ is the subdivision
of three-space by four general planes
through the origin. Corollary \ref{normaltriang}
describes a triangulation where each of the
square-based cones are subdivided into two
triangular cones. For instance, the
normal cone of $\Delta_\FF$ at the vertex $(0,1,2,1)$ 
is subdivided into the normal cones
  of $\Delta_{\widehat \FF}$ at the vertices
$(1,2,7,3)$ and $(1,3,7,2)$.
Likewise, the normal cone of $\Delta_\FF$
at  $(0,2,0,2)$ 
is subdivided into normal cones 
of $\Delta_{\widehat \FF}$ at the vertices
$(1,4,1,7)$ and $(1,7,1,4)$. 
The normal cone  of $\Delta_\FF$ at $(0, 1, 1, 2)$ 
remains unsubdivided. It equals the
normal cone  of $\Delta_{\widehat \FF}$ at $(1, 2, 3, 7)$. \qed
\end{explrm}

We close this section by describing a convenient
representation of nested sets in terms of labeled trees.
This representation appeared implicitly in the
proof of Theorem \ref{einsfuenf}, and in a modified version for partition
lattices it was used previously in~\cite{F}.
Let $T$ be a rooted tree whose nodes $\nu$
are labeled by non-empty pairwise disjoint subsets
$T_\nu$ of $[r]$. For each node $\nu$ we write
$T_{\leq \nu}$ (resp.~$T_{<\nu}$)
for the union of all sets $T_\mu$
where $\mu$ is any node in the subtree of $T$
below the node $\nu$ and including $\nu$ (resp.~excluding $\nu$).
We write $\,{\rm sets}(T)\,$ for the 
set family  $\{T_{\leq \nu}\}_\nu \,$
where $\nu$ is a non-root node of $ T$. Since
$T$ has at most~$r$ nodes, including the root,
the cardinality of ${\rm sets}(T)$ is at most $r-1$,

We now fix a  building set $\FF$ in $2^{[r]}$ with $[r] \in \FF$.
A tree $T$ whose nodes are labeled by non-empty pairwise disjoint subsets
as above is called an {\em $\FF$-tree}
if $\,{\rm sets}(T)\,\subseteq \FF$.

\begin{prop} \label{treerepresentation}
For any nested set $\Ss $ of $\FF$ there exists a unique
$\FF$-tree $T$ such that ${\rm sets}(T) = \Ss$.
The nested set complex equals $\,\NN(2^{[r]},\FF) \,= \, \bigl\{ {\rm sets}(T) \,: \,
\hbox{$T$ is an $\FF$-tree}\,\bigr\}$.
\end{prop}

  \begin{proof}
  If $T$ is an $\FF$-tree then $\,{\rm sets}(T)\,$ is nested because
  $\,T_\nu = T_{\leq \nu} \backslash T_{< \nu}\,$ is non-empty. Conversely,
  suppose that $\Ss$ is nested with $k$ elements and contains $G_1=[r]$.
  We build the tree $T$ inductively, starting with the root.
   Let $G_2,\ldots,G_k$ be the maximal elements of $\Ss \backslash \{[r]\}$.
  We label the root by the non-empty set $\,\rho \, = \,[r] \backslash
   (G_2 \cup \cdots \cup G_k)$. If $\,\Ss=\{[r]\}$, 
then $\,\rho = [r]\,$ and we are done.
   Otherwise note that
    the restriction of $\FF$ to $G_i$ is
   a building set in $2^{G_i}$, and the restriction of the nested set $\Ss$ to $G_i$ is
   a nested set. By induction, it is represented by a labeled tree $T_i$
whose root is labeled by $G_i$. Attaching the labeled trees $T_2,\ldots,T_k$
to the root $\rho$, we obtain the unique $\FF$-tree $T$ with 
${\rm sets}(T) = \Ss$.
  \end{proof}


\section{Triangulations of the Bergman complex} \label{sect_triangBM}

In this section we prove the following theorem
relating nested sets and Bergman fans.

\begin{thm} \label{Main}
For any matroid $M$ and any building set
$\GG$ in its lattice of flats~$\LL_M$,
the fan $\wti \Sigma(\LL_M,\GG)$ refines the Bergman fan $\wti \BB(M) $.
The geometric realization 
$\Sigma(\LL_M,\GG)$ of the nested set complex
$ \NN(\LL_M,\GG) $ is a 
triangulation of the Bergman complex $\BB(M)$.
\end{thm}

We will first prove a local version of Theorem \ref{Main}.
The Bergman complex $\BB(M)$ of a matroid $M$ of rank $r$
is an $(r-2)$-dimensional subcomplex in the boundary of
the dual matroid polytope $P_M^*$. 
We fix a basis $\sigma$ of $M$. The
{\em local Bergman complex} $\BB_\sigma(M)$ is defined as the 
intersection of $\BB(M)$ with the facet
of $P_M^*$ dual to the vertex $e_\sigma$ of 
the matroid polytope $P_M$. Equivalently, 
we can consider the {\em local Bergman fan}
$\wti \BB_\sigma (M)$, which is the
restriction of the Bergman fan $\wti \BB(M)$
to the maximal cone of the normal fan of $P_M$
indexed by $\sigma$.  Consider the sublattice
$\LL_M(\sigma)$ 
of the geometric lattice
consisting of all flats of $M$ that
are spanned by subsets of the basis $\sigma$.
Clearly, $\LL_M(\sigma)$ is a Boolean lattice of rank $r$,
i.e., it is isomorphic to the lattice of subsets
of $\{1,2,\ldots,r \}$.

Let $\GG$ be any building set in $\LL_M$. We write
  $\GG_\sigma$ for the set of all flats in $\GG$ which 
are spanned by subsets of the basis $\sigma$. Then
$\GG_\sigma$ is a building set in the Boolean lattice~$\LL_M(\sigma)$.
The nested set complex 
$\,\NN \bigl(\LL_M(\sigma),\GG_\sigma \bigr)\,$
is called the {\em localization} of
the big nested set complex $\,\NN(\LL,\GG)$ at the
basis $\sigma$. We assume here that the matroid $M$ is connected.

\begin{thm} \label{einseins}
The localization $\,\NN \bigl(\LL_M(\sigma),\GG_\sigma \bigr)\,$
  is a  triangulation of the  local Bergman complex $\BB_\sigma(M)$.
Both complexes are homeomorphic to the $(r-2)$-sphere.
Each of them is naturally realized as the boundary
complex of an $(r-1)$-dimensional polytope.
\end{thm}

\begin{proof} After relabeling we may assume that
the basis $\sigma$ of the matroid $M$ equals
$\sigma = \{1,2,\ldots,r\}$.
Every element $\,i \in \{r+1,r+2,\ldots,n\}\,$
lies in the span of a unique subset
$F_i $ of the basis $\sigma$. This specifies the
following family of subsets of $[r]$:
$$ \FF \quad = \quad \bigl\{ F_{r+1}, F_{r+2}, \ldots, F_n \bigr\}. $$
Let $\,w = (w_1, \ldots,w_r, w_{r+1}, \ldots,w_n )\,$
be a vector in the local Bergman fan
$\wti \BB_\sigma (M)$.  Using the fact that $\sigma = \{1,\ldots,r\}$ is a basis
of the matroid $M_w$, and applying the ``minimum attained twice condition''
to the {\em basic circuit} $\,C = F_i \cup \{i\}$ of $M$, we find that
$$ w_i \,\,= \,\, {\rm min}\{ w_j \,: \, j \in F_i \} \qquad \hbox{for} \,\,\,
  i = r+1,r+2,\ldots,n. $$
This  defines a piecewise-linear map
from $\R^r$ onto the support of the local Bergman fan:
\begin{equation}
\label{mappsi}
  \psi \,: \,\R^r \,\mapsto |\tilde \BB_\sigma (M)|\, ,\,\,
(w_1, \ldots,w_r) \,\mapsto \, 
(w_1, \ldots,w_r, w_{r+1}, \ldots,w_n ). 
\end{equation}
This map is obviously a bijection.
The domains of linearity of the $i$-th coordinate (for $i > r$)
of the map $\psi$ is the normal fan of the simplex
$\Delta_{F_i}$. The common
refinement of these normal fans is the normal
fan of the polytope $\Delta_\FF$. Hence
the domains of linearity of the map $\psi$
are the cones in the normal fan of
$\Delta_\FF$. We conclude that $\psi$ induces
a combinatorial isomorphism between the
normal fan of $\Delta_\FF$ and the local 
Bergman complex $B_\sigma(M)$.

The map $\, 2^{[n]} \rightarrow 2^{[r]} = \LL_M(\sigma), \,F \mapsto [r] \cap F \,$
defines a bijection between flats spanned by
subsets of $\sigma = [r]$ and subsets of $[r]$.
Under this bijection, the flats in $\GG_\sigma$
are identified with subsets of $[r]$, and we have
$$ \FF \, \subseteq  \, \widehat{\FF} \,\subseteq \, \GG_\sigma \,\,\,\subset \,\,\, 2^{[r]}.$$
By Theorem \ref{einsfuenf}, $\,\NN(2^{[r]},\GG_\sigma)\,$
  is the normal fan of the simple polytope $\Delta_{\GG_\sigma}$.

We have shown that both the local Bergman complex $\BB_\sigma(M)$ and
the localization $\,\NN \bigl(\LL_M(\sigma),\GG_\sigma \bigr)\,$
arise as boundary complexes
of $(r-1)$-dimensional polytopes, namely,
the polytopes dual to $\Delta_\FF$ and 
$\Delta_{\GG_\sigma}$, respectively.
Since the polytope $\Delta_\FF$ is a Minkowski summand of $\Delta_{\GG_\sigma}$,
it follows that $\,\NN \bigl(\LL_M(\sigma),\GG_\sigma \bigr)\,$
is a triangulation of $\BB_\sigma(M)$. 
\end{proof}

We are now prepared to prove the theorem stated at the
beginning of this section.

\medskip
\noindent {\sl Proof of Theorem \ref{Main}. }
Both the nested set complex $\,\NN(\LL_M,\GG)\,$
and the Bergman complex $\,\BB(M)\,$
are regarded as polyhedral complexes in the
sphere $\Sphere $, so it suffices to prove the second
assertion. Using the map $\psi$ in (\ref{mappsi}), it is easy to see
that these spherical complexes 
have the same support. Indeed, if $w \in \Sphere$ and
$\sigma$ is any basis of $M_w$, then 
$$ 
w \in |\NN(\LL_M,\GG)|
\,\,\iff \,\,
w \in |\NN(\LL_M(\sigma),\GG_\sigma)|
\,\,\iff \,\,
w \in |\BB_\sigma (M) | 
\,\,\iff \,\,
w \in | \BB(M) |.$$
Since
$\, \NN(\LL_M(\sigma),\GG_\sigma) \,$ triangulates
$\,\BB_\sigma(\GG)$ by Theorem \ref{einseins},
taking the union over all bases $\sigma$ of $M$,
it follows that $\,\NN(\LL_M,\GG)\,$
triangulates $\,\BB(M)\,$.
\qed \medskip

\begin{crl} {\rm (Ardila-Klivans \cite[Sect.\ 2, Thm.\ 1]{AK}) }
The Bergman complex~$\BB(M)$ of a matroid $M$ is homeomorphic to the order
   complex~$\Delta(\LL_M)$ of the lattice of flats~$\LL_M$. In
   particular, $\BB(M)$ is homotopy equivalent to a wedge of spheres of
   dimension~$r{-}2$; the number of spheres is given by the absolute
   value of the M\"obius function on~$\LL_M$.
\end{crl}

\begin{proof}
The first assertion is Theorem \ref{Main}
applied to the largest building set $\GG\,{=}\,\LL$.
That the second assertion holds for the order complex $\Delta(L_M)$
is a well-known result in topological combinatorics (see e.g.~\cite{Bj}).
Hence the first assertion implies the second.
\end{proof}

Theorem \ref{Main} raises the following combinatorial question:
What is the matroid $M_\Ss$ corresponding to a particular
nested set $\Ss$ of $\GG$~? Here $M_\Ss = M_w$, where
$w$ is any point on the sphere $\Sphere$ that lies
in the relative interior of the simplex corresponding to $\Ss$.
To answer this question, we represent the nested set 
$\Ss$ by a labeled tree $T_\Ss$ as follows. 
Fix any basis $\sigma$ such that
$\Ss \in \NN(\LL_M(\sigma), \GG_\sigma)$.
  By Proposition \ref{treerepresentation}, there exists
  a $\GG_\sigma$-tree $T$ with $\Ss = {\rm sets}(T)$.
  Consider the flats
$\, F_{\leq \nu} = {\rm span}(T_{\leq \nu})\,$
and $\, F_{< \nu} = {\rm span}(T_{< \nu}) \,$ of the matroid $M$.
Then $[F_{< \nu}, F_{\leq \nu}]$ is an interval in the geometric 
lattice $\,\LL_M$. We now replace the label $T_\nu$ of the
node $\nu$ in the tree $T$ by this interval. The resulting labeled
tree $T_\Ss$ is independent of the choice of the basis $\sigma$.
It only depends on the nested set~$\Ss$. 
In the following description of the matroid~$M_{\Ss}$ we will denote
the matroid defined by an interval $[F,G]$ in the geometric lattice~$\LL_M$
by $M[F,G]$, extending our notation for restrictions and contractions of 
matroids introduced prior to Proposition~\ref{flacetprop}.

\begin{thm} \label{whatmatroid}
The matroid $M_\Ss$ is the direct sum of the matroids which are
defined by the geometric lattices that appear as labels
of the nodes in the tree $T_\Ss$. In symbols,
\begin{equation}
\label{nomorebullet} M_\Ss \quad = \quad
\bigoplus_{\nu \,{\rm node}\,{\rm of}\,T_\Ss}
M[ F_{< \nu}, F_{\leq \nu}]. 
\end{equation}
\end{thm}

\begin{rem} \label{earlier} \rm
The special case of the formula (\ref{nomorebullet}) where $\Ss = \{F,[r]\}$
represents a vertex of $\,\NN(\LL_M,\GG)\,$ 
appears in equation (\ref{bulletformula})
in the proof of Proposition \ref{flacetprop}.
\end{rem}

\begin{rem} \rm
The case of Theorem \ref{whatmatroid} where
$\GG = \LL_M$ is the maximal building set was
proved by Ardila and Klivans \cite[Section 2]{AK}.
In their work, the tree~$T_\Ss$ is always a chain.
\end{rem}

\medskip 
{\sl Proof of Theorem \ref{whatmatroid}. } A basis $\sigma$ of $M$ is
a basis of $M_\Ss$ if and only if $\Ss \in \NN(\LL_M(\sigma),
\GG_\sigma)$. We can use that basis to construct the tree $T_\Ss$.
The subset $T_\nu$ of $\sigma$ is a basis of the matroid $\,M[F_{<
\nu}, F_{\leq \nu}]\,$ and hence $\sigma$ is a basis of the matroid on
the right hand side of (\ref{nomorebullet}).  Conversely, suppose that
$\sigma$ is a basis of the matroid on the right hand side of
(\ref{nomorebullet}). Then the set $\, {\widetilde T}_\nu\, = \,
(\sigma \cap F_{\leq \nu}) \backslash F_{< \nu}\,$ is a basis of the
matroid $\,M[F_{< \nu}, F_{\leq \nu}]$.  If we take $\widetilde T$ to
be the same tree as $T$ but with each node labeled by ${\widetilde
T}_\nu$, then $\widetilde T$ is a $\GG_\sigma$-tree, and we conclude
that $\,\Ss \in \NN(\LL_M(\sigma), \GG_\sigma)$, i.e.,~$\sigma$ is a
basis of $M_\Ss$.  \qed

\medskip

Of special interest is the case when $\Ss$ is a facet of the nested
set complex $\NN(\LL_M,\GG)$. In that case, $M_\Ss$ is a {\em
transversal matroid}, which means that $M_\Ss$ is a direct sum of
matroids of rank $1$. Indeed, if $\,|\Ss| = r\,$ then $T_\Ss$ is a
binary tree, and the total number of nodes $\nu$ is $r$. For each node
$\nu$, the set $T_\nu$ is a singleton, and the matroid $\,M[F_{< \nu},
F_{\leq \nu}]\,$ has rank~$1$.  Since $\,M[F_{<\nu}, F_{\leq \nu}]$
has no loops, by construction, this rank one matroid is uniquely
specified by the subset $\,F_{\leq \nu} \backslash F_{< \nu}\,$ of
$[n]$.  The collection of sets $\,F_{\leq \nu} \backslash F_{< \nu}$, 
as $\nu$ runs over the nodes of $T_\Ss$, is a partition of
the set $[n]$ in $r$ parts. The
transversal matroid $M_\Ss$ is uniquely specified by this set partition.

\begin{crl} \label{crl_facetsN}
The facets of the nested set complex $\NN(\LL_M,\GG)$
are indexed by partitions of the ground set $[n]$ into $r$ parts.
\end{crl}


\section{Computing nested sets and the Bergman complex} \label{sect_NMvBM}

In this section we consider a connected matroid $M$ of rank $r$ on
$n$ elements,
and we fix $\,\GG = \GG_{\rm min} \,$ to be the {\em minimal\/} 
building set in its lattice of flats~$\LL_M$. We denote the corresponding 
nested set complex by
\[ 
\NN(M) \quad = \quad \NN( \LL_M, \GG_{\rm min})\, , 
\] 
and, for the
purpose of this section, we call $\NN(M)$ simply {\bf the} {\em nested
set complex} of  $M$.

The minimal building set $\GG_{\rm min}$ consists of all connected
flats $F$ of $M$. Observe that this description coincides with the set
of {\em connected\/} elements in $\LL_M$ that we identified as the minimal
building set for arbitrary lattices in the beginning of
Section~\ref{sect_ns}.

By our results in Section 4, the various nested set complexes  triangulate 
the Bergman complex, and by Proposition~\ref{prop_subd}, the 
minimal nested set complex~$\NN(M)$ is the coarsest among these 
triangulations. 
The vertices of the Bergman complex $\BB(M)$ are the
flacets of the matroid (by Theorem \ref{dfbergman}), and the vertices of
the nested set complex $\NN(M)$
are the connected flats of $M$ (by definition). These notions do not 
coincide in general:

\begin{rem} \label{rm_flconnfl}
  \rm Every flacet is a connected flat but not conversely.
The matroid $M'$ in Examples \ref{K4minus}
and \ref{ex_NK4minus}  has a connected
flat of rank one (i.e.,~a single element) which is not a flacet:
the edge which is complementary to the edge 
that was removed from $K_4$.
\end{rem}

This distinction between flacets and connected flats amounts to the
fact that in general new vertices are added when passing from the
Bergman to the nested set complex:

\begin{crl} The Bergman complex $\BB(M)$ is triangulated
by the nested set complex $\NN(M)$ without 
additional vertices if and only if  every connected
flat of $M$ is a flacet.
\end{crl}

We shall now generalize this corollary to give
a criterion for when $\NN(M)$ equals $\BB(M)$.

\begin{thm} \label{whenequal} The nested set complex $\NN(M)$ equals
  the Bergman complex $\BB(M)$ if and only if
the matroid $M[F,G]$ is connected for
  every pair of flats $F \subset G$ with $G$  connected.
\end{thm}

\begin{proof}
Consider any simplex $\Ss$ of the nested set complex $\NN(M)$
and let $\Gamma_\Ss$ be the smallest face of the
Bergman complex $\BB(M)$ which contains $\Ss$.
A necessary and sufficient condition for $\NN(M) = \BB(M)$ 
to hold is that ${\rm dim}(\Ss) = {\rm dim}(\Gamma_\Ss)$ for
all such pairs $\Ss \subseteq \Gamma_\Ss$.
Now, ${\rm dim}(\Ss)$ is simply $|\Ss| - 1$, and
$\,{\rm dim}(\Gamma_\Ss) \,$ equals $\, c(M_\Ss)-1$. From this we conclude
$$ \centerline{$\NN(M) = \BB(M)\,\,\,$ if and only if
$\,\,\,\, c(M_\Ss) = |\Ss| \,$ for all $\Ss \in \NN(M)$.} $$

The matroid $M_\Ss$ was characterized in
Theorem \ref{whatmatroid}. We have
$\,c(M_\Ss) = |\Ss| \,$ if and only if
all the matroids $\,M[ F_{< \nu}, F_{\leq \nu}]\,$
in the decomposition (\ref{nomorebullet})
are connected. Note that here $F_{\leq \nu}$
is always a connected flat and $F_{< \nu}$ 
is a subflat of $F_{\leq \nu}$. This establishes
the if-direction of Theorem \ref{whenequal}.
For the only-if direction, we consider any
  pair of flats $F \subset G$ such that $G$ is connected.
  Let $\Ss$ denote the nested set which consists of $G$
  and the connected components of $F$. If $\nu$ is the
  root of the tree $T_\Ss$ then we have
$\,F_{< \nu} = F\,$ and  $\, F_{\leq \nu} = G$.
This shows that  $\NN(M) = \BB(M)$ implies
the connectedness of $M[F,G]$.
\end{proof}

\begin{rem} \label{SpaceOfTrees} \rm
Consider the graphic matroid $M = M(K_n)$ whose bases
are the spanning trees in the complete graph $K_n$.
Here the criterion of Theorem \ref{whenequal}
is satisfied, and the nested set complex
coincides with the Bergman complex. This was seen 
for $n=4$ in Examples~\ref{K4} and \ref{ex_NPi4Gmin},
and in~\cite[Rem.\ 3.4.(2)]{F} for general $n$.
Ardila and Klivans \cite[Sect.~3, Prop.]{AK} showed that
  $\,\wti \NN(M) = \wti \BB(M)\,$ equals the 
{\em space of phylogenetic trees}.
Theorem~\ref{whatmatroid} states that
{\bf every} nested set complex
can be interpreted as a certain complex of trees.

The hyperplane arrangement  corresponding to $M$
is the {\em braid arrangement}
$\{x_i = x_j\}$. What we are discussing here
is its {\em wonderful model} (in
the sense of De~Concini and Procesi~\cite{DP},
see Section 6). Ardila, Reiner and Williams~\cite{ARW} recently
showed that nested set complexes and Bergman complexes coincide 
for any finite reflection arrangement.
It might be interesting to classify all subarrangements
of reflection arrangements (e.g.,~graphic matroids)
for which the Bergman complex equals the nested set complex.
\end{rem}

\smallskip

We next present an algorithm for
computing both the Bergman complex $\BB(M)$
and its triangulation by the nested set complex $ \NN(M)$.
We prepared a test implementation of this algorithm 
in  {\tt maple}. This
code is available from Bernd Sturmfels upon request.

\begin{alg} \label{maplealgo}
  (Computing the Bergman complex and the nested set complex). 
\hfill \break \rm
Input: A rank $r$ matroid $M$ on $[n]$, given by its bases. \hfill \break 
Output: All maximal faces of the Bergman complex $\BB(M)$,
represented by unordered partitions $\{B_1,\ldots,B_r\}$ of $[n]$ into $r$ 
non-empty blocks, as described in Corollary~\ref{crl_facetsN}.

\smallskip

\noindent
1. Initialize $\Omega = \emptyset$.
  \phantom{dadadadad} 
   {\it (This will later be the set of all unordered partitions)}.

\noindent 2. Precompute all connected flats of $M$.

\noindent
3. For every basis $\sigma$ of $M$ do the following:
\begin{itemize}
\item[(a)]
For each $\,i \in [n] \backslash \sigma \,$ find the unique
set $\,F_i \subseteq \sigma\,$ such that $\,F_i \cup \{i\}\,$ is a circuit. \break
\phantom{dadadadadadadadada} {\it (The local building set equals 
$\,\GG_\sigma \,=\,
\bigl\{\,F_i \,:\, i \in [n] \backslash \sigma \,\bigr\}$).}
\item[(b)]
For each permutation $\pi$ of $[r]$ do
\begin{itemize}
\item[(i)]
for each $j \in [r]$ set
$\,\omega_j \, := \, \{j \}\,\cup \,
\{\,i \in [n]\backslash \sigma  \,: \,{\rm min}(F_i) = \pi_j \,\} $.
\item[(ii)] Set $\,\Omega \,:= \,\Omega \,\cup \,
\bigl\{ \{\omega_1,\omega_2 ,\ldots,\omega_r \} \bigr\}$.
\end{itemize}
\end{itemize} 
\noindent 4. For each $\, \omega \in \Omega \,$ do:
\begin{itemize}
\item[(a)] Output: ``The partition $\,\omega \,$ represents a facet of $\BB(M)$''.
\item[(b)] Set $\,\Pi := \emptyset$. \phantom{dadada}
{\it (This will be the nested set triangulation of the facet $\omega$)}.
\item[(c)] For each connected flat $F$ of $M$ do
  \begin{itemize}
  \item[$\bullet$] Set $\,s = {\rm rank}(F)$.
  \item[$\bullet$]
  If $F$ is the union of  $s$ blocks $\,\omega_{i_1},
  \ldots,\omega_{i_s}\,$ then $\,\Pi \, := \,\Pi \,\cup \,\{F\}$.
  \end{itemize}
  \item[(d)] If the cardinality of $\Pi$ is $r$ then  output \hfill \break
  ``The simplex $\omega$ is not subdivided; it equals the simplex $\Pi$ in $\NN(M)$.''
  \item[(e)] Otherwise compute and output the set of all maximal nested  sets
  on $\Pi$.
  \end{itemize}

  \end{alg}

\smallskip

\noindent {\sl Discussion and Correctness. }
The loop in Step 3 computes the local Bergman complex~$\BB_\sigma(M)$
for every basis $\sigma$ of $M$, and it saves
all set partitions representing facets of~$\BB(M)$
(as in Corollary~\ref{crl_facetsN}) in one big  set $\Omega$.
Step (b) takes advantage of the many-to-one correspondence between 
permutations $\pi$ of $[r]$ and vertices of $\Delta_{\GG_\sigma}$,
coming from the fact that $\Delta_{\GG_\sigma}$ is a Minkowski
summand of the permutohedron.
In Step 4 we output each facet $\omega$ of $\BB(M)$
along with the list of maximal simplices in its
nested set triangulation. The crucial 
step is the second bullet $\bullet$ in
  Step 4 (c), which tests, for each
connected flat $F$ of $M$, whether or not
the corresponding vertex of $\NN(M)$ lies
on the facet $\omega$ of $\BB(M)$.
\qed

\smallskip

\begin{rem} \rm
We believe that an improved version of
Algorithm  \ref{maplealgo}  should be able to
list the facets of both $\BB(M)$ and $\NN(M)$
in a {\em shelling order}, so as to reveal
the topology of these spaces and to offer
a practical tool for the computation of residues
using the methods of \cite{DP2}.
The  constructing of such shellings is the
goal of a subsequent project.
\end{rem}

\medskip

We ran Algorithm \ref{maplealgo} on
a range of matroids of various ranks, including the 
following two examples which can serve as test cases
for future  tropical algebraic geometry software.

\begin{explrm} \label{ex_R10}
We consider the famous self-dual unimodular matroid 
$R_{10}$ of rank $r\,{=}\,5$ on $n\,{=}\,10$ elements. This matroid plays
a special role in Seymour's decomposition theory for
regular matroids (see \cite{Tru}). The ground set for $R_{10}$
is the set of the edges of the complete graph $K_5$, and its
circuits are the four-cycles and their complements. Its geometric
lattice has $45$ lines, $75$ planes and $30$ three-dimensional subspaces.
There are $40$ connected flats:
\begin{itemize}
\item[$\circ$] the ten points in the ground set,
\item[$\circ$] the 15 four-cycles  (these are planes),
\item[$\circ$] the 5 copies of the complete graph $K_4$ (these are three-spaces),
\item[$\circ$] the 10 copies of the complete bipartite $K_{2,3}$ (also three-spaces).
\end{itemize}

Topologically, 
the Bergman complex of $R_{10}$ is a bouquet of nine $3$-dimensional spheres.
It is constructed as follows. We first note that the nested set complex 
of $R_{10}$ consists of $405$ tetrahedra, which come in eight families:

Each of the five copies of $K_4$ contributes 27=1+18+4+4 tetrahedra:

\begin{itemize}
\item[(a)]  the four edges not in $K_4$,
\item[(b)]  the $K_4$, a four-cycle in $K_4$ and any two of its edges,
\item[(c)]  the $K_4$, and three of its edges that form a $K_3$,
\item[(d)]  the $K_4$, and three of its edges that form a $K_{1,3}$.
\end{itemize}

Each of the ten copies of $K_{2,3}$ contributes 27=1+18+2+6 tetrahedra:
\begin{itemize}
\item[(e)]  the four edges not in $K_{2,3}$
\item[(f)]  the $K_{2,3}$, a four-cycle in $K_4$ and any two of its edges
\item[(g)]  the $K_{2,3}$, and three of its edges that form a $K_{1,3}$
\item[(h)]  the $K_{2,3}$, and three edges that touch all six vertices.
\end{itemize}

The Bergman complex of $R_{10}$ has $360$ facets, namely
$315$ tetrahedra and $45$ bipyramids. The are $15$ bipyramids
formed by pairs of tetrahedra of type (b), namely, two disjoint 
edges in a $K_4$ and the two four-cycles of $K_4$ containing them.
The other $30$ bipyramids are formed by pairs of tetrahedra of type (f),
namely, two edges in $K_{2,3}$ complementary to a four-cycle and the 
two other four-cycles of in $K_{2,3}$  which contain them.
\end{explrm}

\begin{explrm} \label{cographic}
Let $M$ be the {\em cographic matroid} $M(K_5)^*$.
Here $r = 6$ and $n = 10$. The bases of $M$
are the six-tuples of edges in the complete graph $K_5$
which are complementary to the spanning trees.
The nested set complex $\NN(M)$ is a four-dimensional
simplicial complex with $f$-vector $(25, 185, 615, 955, 552)$.
The Bergman complex $\BB(M)$ has $447$ facets,
of which $105$ are subdivided into pairs of
$4$-simplices when passing to $\NN(M)$.
\end {explrm}

\begin{explrm} \label{cubes}
The equality $\BB(M) = \NN(M)$ holds 
when $M$ is the rank $r=4$ matroid specified by
the $n=8$ vertices of the three-dimensional
unit cube. The simplicial complex $\BB(M) = \NN(M)$
has $20$ vertices, $76$ edges and $80$ triangles.
On the other hand, for the four-dimensional
unit cube $(r=5,n=16)$, the Bergman complex $\BB(M)$ is not
simplicial (it has $2600$ facets). It is properly subdivided
by the nested set complex $\NN(M)$ which consists of
  $176 $ vertices, $1280 $ edges, $3360 $ triangles and
$2720 $ tetrahedra.
\end {explrm}


\section{Tropical compactification of arrangement complements}

In this section we interpret our combinatorial results in terms of
algebraic geometry. We consider a connected rank $r$ matroid $M$ on
$[n]$, so its matroid polytope $P_M\,{\subset}\,\R^n$ has dimension $n{-}1$,
and we assume that $M$ is realized by an $r$-dimensional linear
subspace $V$ of the vector space $\C^n$, or equivalently, by a
$(r{-}1)$-dimensional projective linear subspace $\overline{X}$ in the
complex projective space $\P^{n-1}$.  A subset $D$ of $[n]$ is
dependent in the matroid~$M$ if and only if there exists a non-zero
linear form $\,\sum_{i \in D} c_i x_i\,$ that vanishes on $\overline{X}$.

We identify the algebraic 
torus $(\C^*)^{n-1}$ with the projective space $\P^{n-1}$
minus its $n$ coordinate hyperplane $\,\{x_i=0\}$.
We are interested in the non-compact variety 
$\,X \, = \, \overline{X} \,\cap \, (\C^*)^{n-1}$. This is the complement
of an arrangement of $n$ hyperplanes in $\overline{X} \simeq \P^{r-1}$.

A standard problem in algebraic geometry is to
construct a smooth compactification of the arrangement
complement $X$ which has better
properties than the ambient projective space $\overline{X}$.
Ideally, one wants the complement of $X$ in that
compactification to be a {\em normal crossing divisor}.
A solution to this problem was given by
De~Concini and Procesi \cite{DP}. 

We briefly describe the construction of their {\em wonderful
  compactification} $X_{wond}$. The geometric lattice $\LL_M$ of the
matroid $M$ is the {\em intersection lattice} of the hyperplane
arrangement. Each flat $F \subset [n]$ of rank $s$ in $M$ corresponds
to a subspace \mbox{$\{ x_i = 0 \, : \, i \in F\}$} of codimension $s$
in the arrangement.  Let $\P^F = \P(\C^F)\,$ denote the coordinate
subspace of $\P^{n-1}$ with coordinates $\,x_i, \, i \in F$, and
consider the projection $\,\psi_F \, : \, \P^{n-1} \rightarrow \,\P^F$.
The restriction of $\psi_F$ to $X$ is a regular map. Let $\GG =
\GG_{min}$ be the minimal building set which consists of the connected
flats, and consider the product of all of these regular maps
\begin{equation}
\label{DPmap} \psi_\GG \,\, : \,\, X \,\rightarrow \, 
\prod_{G \in \GG} \P^G\,, \quad
u \, \mapsto \, \bigl(\psi_G(u)\,:\,G \in \GG \bigr). 
\end{equation}
The {\em wonderful model} $X_{wond}$ is the closure of $\,\psi_\GG(X)\,$  in
the compact variety $\,\prod_{G \in \GG} \P^G$.

More recently, Tevelev \cite{Tev} introduced an alternative
compactification, called the {\em tropical compactification}
$\,X_{trop}$. The advantage of this new construction is that $X$ can
now be an arbitrary subvariety of the algebraic torus $(\C^*)^{n-1}$.
Let $\overline{X}$ be the closure of $X$ in~$\P^{n-1}$ and let ${\rm
  Hilb}_X(\P^{n-1})$ be the Hilbert scheme of all subschemes of $\P^n$
which have the same Hilbert polynomial as $\overline{X}$. Without loss
of generality, we may assume that $X$ is not fixed by any non-unit
element $t \in (\C^*)^{n-1}$. Then the following map is an embedding:
\begin{equation}
\label{isanemb} \phi \,: \,X \rightarrow {\rm Hilb}_X(\P^{n-1}) , \,\,\, t \,\mapsto
\,t^{-1} \cdot X. 
\end{equation}
The tropical compactification $X_{trop}$ is the closure of $\phi(X)$
in the compact scheme $\,{\rm Hilb}_X(\P^{n-1})$. The combinatorial
structure of $X_{trop}$ is governed by the {\em tropical variety}
of~$X$. This is a subfan of the Gr\"obner fan of the ideal $I_X$ of
$X$. This subfan consists of all cones such that the corresponding
initial ideal contains no monomials. These cones correspond to the
strata in the boundary $\,X_{trop} \backslash X$. For nice varieties
$X$, the divisor $\, X_{trop}\backslash X$ is much better behaved than
$\,\overline{X} \backslash X$, and it is expected that $X_{trop}$ is
smooth and its boundary $X_{trop}\backslash X$ can be made normal
crossing using only toric blow-ups \cite{Hak, Tev}.

Here we consider the nicest case, when $\overline{X} \subset
\P^{n-1}$ is a linear space (with matroid~$M$) 
and $X$ is the complement of the
arrangement of the  $n$ hyperplanes  $\{x_i=0\}$ in $\overline{X}$.
The subspaces of $\P^{n-1}$ that can  be gotten by intersecting a 
subset of the $n$ hyperplanes correspond to the proper flats $F$
of the matroid $M$.  Note that a flat $F$ of $M$ is connected
 if and only if it
corresponds to a divisor in the wonderful
compactification $\,X_{wond}$.

\begin{thm} \label{wonderfultotrop}
For any hyperplane arrangement complement $X$,
there is a canonical morphism  from
the wonderful compactification $\,X_{wond}$
onto the tropical compactification $X_{trop}$.
This morphism is an isomorphism whenever the following
combinatorial condition holds:
If $G \in \LL_M$ corresponds to a divisor in $X_{wond}$
then $G$ also corresponds to a divisor in
$\,(X\,\cap \,F)_{wond}\,$
where $F$ is any intersection of hyperplanes  which contains $G$.
\end{thm}

The combinatorial condition above is a translation of
the condition in Theorem \ref{whenequal}. However,
as Remark~\ref{keelremark} shows, we must apply
the condition in Theorem~\ref{wonderfultotrop}
to the cone over the  arrangement, because
$G$ might be the full set $[n]$ in Theorem~\ref{whenequal}

\begin{expl} \rm
Consider an arrangement of five planes 
in $\P^3$ which intersect in one point $G$
and which represent the matroid $M'$ in Examples \ref{K4minus} 
and \ref{ex_NK4minus}. Let $F$ be the plane indexed by
the edge $5$ in Figure \ref{fig_NK4minus}. Then $G$
does not correspond to a divisor in  $\,(X\,\cap \,F)_{wond}\,$
because the restriction of the arrangement to $F$
consists of two lines through~$G$.
Hence $\, X_{wond} \rightarrow X_{trop}\,$ is not
an isomorphism.
On the other hand, the recent work of
Ardila, Reiner and Williams \cite{ARW}
implies that, for any finite reflection arrangement,
the tropical compactification coincides with the
wonderful compactification.
\end{expl}

\begin{rem} \label{keelremark} \rm
  The case when $\overline{X}$ is a two-dimensional plane in
  $\P^{n-1}$ is completely described in the work of Keel and Tevelev
  \cite[Lemma 8.11, p.~43]{KT}. Here we are considering an arrangement
  of $n$ lines in $\P^2$.  The wonderful compactification is equal to the
   tropical compactification except when there is a
  line $L$ in the arrangement and two points $a,b\in L$ such that each
  remaining line passes through $a$ or $b$.  In that case,
  $X_{trop}$ is obtained by blowing up the points $a,b \in \P^2$, then
  contracting the strict transform of the line $L$. The tropical
  compactification is therefore isomorphic to $\P^1 \times \P^1$, and
  the strict transforms of the lines through $a$ and $b$ become fibers
  of the first and second projections, respectively.
\end{rem}

In what follows we shall give a more concrete description
of the tropical compactification $X_{trop}$ and we shall
explain how the morphism $\,X_{wond} \rightarrow X_{trop}\,$ works.
Since $V$ is a  linear space, the
Hilbert scheme ${\rm Hilb}_X(\P^{n-1})$ is simply the
Grassmannian $Gr_{r,n}$ of $r$-dimensional
linear subspaces in $\C^n$.  Let $I_X$ denote the ideal
in $\C[x_0,x_1,\ldots,x_n]$ generated by all linear
forms that vanish on $X$. For  $t = (t_1: \cdots : t_n)  \in (\C^*)^{n-1}$,
we write $\,t \cdot I_X\,$ for the image of $I_X$ under replacing
$x_i$ by $t_i x_i$ for all $i$. Clearly, $\,I_{t ^{-1}\cdot X} \,=\,
t \cdot I_X \,$ for any $t \in (\C^*)^{n-1}$, i.e.,
the variety of the ideal $\,t \cdot I_X\,$ is the
  translated subspace $\, t^{-1} \cdot X $.

The map $\, t\, \mapsto \, t \cdot I_X \,$ defines an embedding
of the torus $(\C^*)^{n-1}$ into the Grassmannian  ${\rm Gr}_{r,n}$.
The closure of its image is a projective toric variety $T_M$.
As an abstract toric variety $T_M$ depends only on the matroid $M$,
namely, it is the toric variety associated with the matroid polytope $P_M$.
However, the specific embedding of $T_M$ into ${\rm Gr}_{r,n}$
depends on the specific realization $X$ of $M$.
The coordinate ring of $T_M$ is the {\em basis monomial ring} \cite{Wh0}
of the matroid $M$, which is the subalgebra of $\,\C[t_1,\ldots,t_n]\,$
generated by the monomials $\,\prod_{i \in \sigma} t_i\,$ where
$\sigma$ runs over all bases of $M$.
A  result of Neil White \cite{Wh0} states that the basis monomial ring
is a normal domain, hence the toric variety $T_M$ is arithmetically normal.

The $(\C^*)^{n-1}$-action on the Grassmannian ${\rm Gr}_{r,n}$ restricts to
the action of the dense torus on the toric variety $T_M$.
(Note that ${\rm dim}(T_M) = n-1$ because $M$ is connected). 
The  $(\C^*)^{n-1}$-orbits on $T_M$  correspond to the distinct
{\em initial ideals} of the  ideal $I_X$, i.e.,
\begin{equation}
\label{InitialIdeal}
  {\rm in}_w(I_X) \,\,\, = \,\,\, \langle \, {\rm in}_w(f) \,: \, f \in I_X
\,\rangle .
\end{equation}
Two vectors $w, w' \in \R^{n}$ lie in the same cone of
the fan of the toric variety $T_M$ if and only if $\, {\rm in}_w(I_X)
= {\rm in}_{w'}(I_X)$. This fan is the normal fan of the matroid polytope $P_M$.

We now restrict the map  $\, t\, \mapsto \, t \cdot I_X \,$ 
from $(\C^*)^{n-1}$ to its subvariety $X$. The result is
precisely the map $\phi$ defined in (\ref{isanemb})
above. We can now rewrite this map as follows:
\begin{equation}
\label{phimap} \, \phi\, :\, X \ \rightarrow \, T_M \,, \,\,\, t \,\mapsto \, t \cdot I_X . 
\end{equation} \nopagebreak[4] 

\begin{rem} \rm
The map $\,\phi\,$ defines an embedding of the 
hyperplane arrangement complement $X$ into the
projective toric variety $T_M$ associated with  the matroid polytope $P_M$.
The tropical compactification $X_{trop}$ is the closure of the
image of this embedding.
\end{rem}

The following lemma characterizes the location of $X_{trop}$ inside $T_M$
combinatorially.

\begin{lm}
Let $C$ be a cone in the normal fan of the matroid polytope $P_M$,
and let $w $ be any vector in the relative interior of $C$.
Then the following are equivalent:
\begin{enumerate}
\item $C$ is a cone in the Bergman fan $\wti \BB(M)$.
\item The initial ideal $in_w(I_X)$ in (\ref{InitialIdeal})
contains no monomial.
\item The initial ideal $in_w(I_X)$ in (\ref{InitialIdeal})
contains no variable $x_i$.
\item $\,X_{trop}\,$ does not intersect the
$(\C^*)^{n-1}$-orbit of $T_M$ corresponding to $C$.
\end{enumerate}
\end{lm}

\begin{proof}
Statements (2) and (3) are equivalent because
$I_X$ is generated by linear forms. The matroid
associated with the linear ideal ${\rm in}_w(I_X)$
is the matroid $M_w$, and (3) means that the
element $i$ is not a loop of $M_w$. Thus
(3) is equivalent to (1). The equivalence of 
(2) and (4) holds in general for tropical
compactifications of arbitrary varieties $X \subset (\C^*)^n$.
It is seen by considering the universal family
over the Hilbert scheme in (\ref{isanemb}).
\end{proof}


\noindent {\sl Proof of Theorem \ref{wonderfultotrop}. }
Every projective space $\P^G$ in 
(\ref{DPmap}) corresponds to a face $\Delta_G$
of the $(n-1)$-dimensional simplex  $\Delta$.
Let $\GG$ be the collection of all connected flats
and let $\, \Delta_\GG\,$ denote the Minkowski sum
of the simplices $\Delta_G$, where $\,G \in \GG$.
Let $T_\GG$ denote the projective toric variety
associated with the polytope $\Delta_{\GG}$.
The toric variety~$T_\GG$ is gotten from projective
space $\,\P^{n-1}$ by the sequence of blow-ups
along the linear coordinate subspaces indexed by $\GG$.
The map $\psi_\GG$ in (\ref{DPmap}) can therefore be replaced
by the inclusion $\, X \subset (\C^*)^{n-1} \subset T_\GG$,
i.e., the wonderful compactification $X_{wond}$ coincides with
  the closure of the arrangement complement $X$ in the
  projective toric variety $T_{\GG}$. Note that $T_\GG$ is
  generally not smooth, but $X_{wond}$ only meets
  smooth strata of $T_\GG$, which explains all
  the  wonderful properties of this compactification.
  
  Our key observation in this proof is that the matroid polytope $P_M$
  is a Minkowski summand of the polytope $\Delta_\GG$ constructed in
  the previous paragraph. Indeed, $P_M$ is given by the inequalities
  in Proposition \ref{matroidineqs}, where $F$ runs over $\GG$, and
  $\Delta_\GG$ is defined the inequalities in Proposition
  \ref{ineqrepp} where $G$ runs over $\GG$. A careful examination of
  these inequality representations reveals that $P_M$ is a Minkowski
  summand of $\Delta_\GG$.

  The relation ``is Minkowski summand of'' among lattice polytopes
  correspond to projective morphisms in toric geometry.  Since $P_M$
  is a Minkowski summand of $\Delta_\GG$, we get a projective morphism
  from the toric variety $\,T_\GG \,$ onto the toric variety $\,T_M$.
  The inclusion of $X$ in $T_M$ given in (\ref{phimap}) coincides with
  the composed map $\, X \subset (\C^*)^{n-1} \subset T_\GG \rightarrow T_M $.
  Hence the closure of the image of $X$ in $T_\GG$ is mapped by a
  projective morphism onto the closure of the image of $X$ in $T_M$.
  This is the desired morphism $\,X_{wond} \rightarrow X_{trop}$.

Consider the collection of cones $C$ in the normal fan of $\Delta_\GG$
such that $X_{wond}$ intersects the orbit indexed by $C$. This subfan
of the normal fan is precisely the geometric realization of the nested
set complex $\,\NN(M)$. Likewise, the collection of cones $C$ in the
normal fan of $P_M$ such that $X_{trop}$ intersects the orbit
indexed by $C$ is precisely the Bergman fan $\BB(M)$.
If these fans are equal then the projective
morphism $\,X_{wond} \rightarrow X_{trop}\,$
is an isomorphism. Hence the second assertion of
  Theorem \ref{wonderfultotrop} follows from Theorem \ref{whenequal}.
\qed

\bigskip

\noindent {\bf Acknowledgements: }
We are grateful to Paul Hacking for his helpful comments on a
draft of this paper, and to 
Lauren Williams and Alex Postnikov for pointing us
to the references \cite{ARW, CD, Pos}.
The IAS/Park City Mathematics Institute (PCMI, July 2004) and 
the Mathematical Sciences Institute in Berkeley (MSRI) 
provided the setting for us to work on this project.
We are grateful to both institutions for their support.
Bernd Sturmfels was supported in part by the 
U.S.~National Science Foundation (DMS-0200729).

\end{document}